\documentclass[aos,citesort,seceqn,dvips]{arximspdf}
\usepackage{graphics}
\doi{10.1214/09-AOS768}
\volume{38}
\issue{3}
\pubyear{2010}
\firstpage{1845}
\lastpage{1884}

\makeatletter

\newproclaim{definition}{Definition}[section]

\newtheorem{lemma}{Lemma}[section]

\newproclaim{example}{Example}[section]

\newtheorem{proposition}{Proposition}[section]

\newproclaim{assumption}{Assumption}[section]

\newtheorem{theorem}{Theorem}[section]
\newtheorem{corollary}{Corollary}[section]

\newcommand{\dint}{\int\!\!\int}

\makeatother

\begin{document}
\begin{frontmatter}

\title{Weakly dependent functional data}
\runtitle{Dependent functional data}

\begin{aug}
\author[A]{\fnms{Siegfried} \snm{H\"ormann}\corref{}\ead[label=e1]{shormann@ulb.ac.be}} and
\author[B]{\fnms{Piotr} \snm{Kokoszka}\thanksref{t2}\ead[label=e2]{piotr.kokoszka@usu.edu}}
\runauthor{S. H\"ormann and P. Kokoszka}
\affiliation{Universit\'e Libre de Bruxelles and Utah State University}
\address[A]{Universit\'e Libre de Bruxelles\\
CP 210, local O.9.115\\
Bd du Triomphe\\
B-1050 Bruxelles\\
Belgium\\
\printead{e1}}
\address[B]{Department Mathematics and Statistics\\
Utah State University\\
3900 Old Main Hill\\
Logan, Utah 84322-3900\\
USA\\
\printead{e2}}
\end{aug}

\thankstext{t2}{Supported in part by
NSF Grants DMS-08-04165 and DMS-09-31948 at Utah State University.}

\pdfauthor{Siegfried Hormann, Piotr Kokoszka}

% HISTORY:
\received{\smonth{6} \syear{2009}}
\revised{\smonth{11} \syear{2009}}

% ABSTRACT
%
\begin{abstract}
Functional data often arise from measurements on fine time grids
and are obtained by separating an almost continuous time record into
natural consecutive intervals, for example, days. The functions thus
obtained form a functional time series, and the central issue in the
analysis of such data consists in taking into account the temporal
dependence of these functional observations. Examples include daily
curves of financial transaction data and daily patterns of
geophysical and environmental data. For scalar and vector valued stochastic
processes, a large number of dependence notions have been proposed,
mostly involving mixing type distances between $\sigma$-algebras. In
time series analysis, measures of dependence based on moments have
proven most useful (autocovariances and cumulants). We introduce a
moment-based notion of dependence for functional time series which
involves $m$-dependence. We show that it is applicable
to linear as well as nonlinear functional time series. Then we
investigate the impact of dependence thus quantified on several
important statistical procedures for functional data. We study the
estimation of the functional principal components, the long-run
covariance matrix, change point detection and the functional linear
model. We explain when temporal dependence affects the results obtained for
i.i.d. functional observations and when these results are robust to
weak dependence.
\end{abstract}

% KEYWORDS
%

\end{frontmatter}

%s1 ###
\section{Introduction} \label{s:sie}
The assumption of independence is often too strong to be realistic
in many applications, especially if data are collected sequentially
over time. It is then natural to expect that the current observation
depends to some degree on the previous observations.
This remains true for functional data and has motivated the
development of appropriate functional time series models.
The most popular model is the autoregressive model of
Bosq \cite{bosq2000}. This model and its various extensions are
particularly useful for prediction (see, e.g.,
Besse, Cardot and Stephenson \cite{bessecardotstephenson2000}
Damon and Guillas \cite{damonguillas2002},
Antoniadis and Sapatinas \cite{antoniadissapatinas2003}).
For many functional time series it is, however, not clear what
specific model they follow, and for many statistical
procedures it is not necessary to assume a specific model.
In such cases, it is important to know what the effect of
the dependence on a given procedure is. Is it robust to temporal
dependence, or does this type of dependence introduce a serious bias?
To answer questions of this type, it is essential to quantify the
notion of temporal dependence. For scalar and vector time series,
this question has been approached from a number of angles, but,
except for the linear model of Bosq \cite{bosq2000}, for functional
time series data no general framework is available. Our goal in this
paper is to propose such a framework, which applies to both linear
and nonlinear dependence, develop the requisite theory and apply it
to selected problems in the analysis of functional time series.
Our examples are chosen to show that some statistical procedures
for functional data are robust to temporal dependence as quantified
in this paper, while other require modifications that take this dependence
into account.

%f1 ###
%
\begin{figure}[b]

\includegraphics{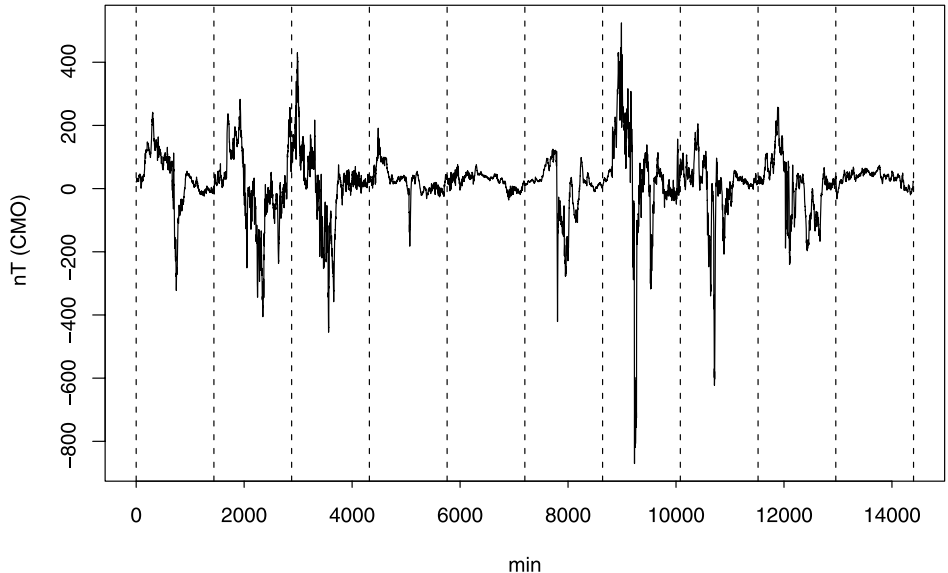}

\caption{Ten consecutive functional observations of
a component of the magnetic field recorded at College, Alaska.
The vertical lines separate days. Long negative spikes lasting
a few hours correspond to the aurora borealis.}
\label{fig:cmo}
\end{figure}

While we focus here on a general theoretical framework, this research
has been motivated by our work with functional data
arising in space physics and environmental science.
For such data, especially for the space physics data, no
validated time series models are currently available, so to justify
any inference drawn from them, they must fit into a general, one might
say, nonparametric,
dependence scheme. An example of space physics data is shown in
Figure \ref{fig:cmo}. Temporal dependence from day to day can be
discerned, but has not been modeled.

The paper is organized as follows. In Section \ref{ss:wd} we introduce
our dependence condition and illustrate it with several examples. In
particular, we show that the linear functional processes fall into our
framework, and present some nonlinear models that also do. It is now
recognized that the \textit{functional principal components} (FPCs) play
a far greater role than their multivariate counterparts (Yao and Lee
\cite{yaolee2006}, Hall and Hosseini-Nasab \cite{hallhn2006},
Reiss and Ogden \cite{reissogden2007}, Benko, H{\"a}rdle and Kneip
\cite{benkohardlekneip2009},
M{\"u}ller and Yao \cite{mulleryao2008}). To develop theoretical
justification for
procedures involving the FPCs, it is necessary to use the convergence
of the estimated FPCs to their population counterparts. Results of
this type are available only for independent observations
(Dauxois, Pousse and Romain \cite{dauxois1982}, and linear processes,
Bosq \cite{bosq2000},
Bosq and Blanke \cite{bosqblanke2007}). We show in Section \ref
{ss:conv} how the
consistency of the estimators for the eigenvalues and eigenfunctions
of the covariance operator extends to dependent functional data.
Next, in Section \ref{s:lr}, we turn to the estimation of an
appropriately defined long-run variance matrix for functional data.
For most time series procedures, the long-run variance plays a role
analogous to the variance--covariance matrix for independent
observations. Its estimation is therefore of fundamental importance,
and has been a subject of research for many decades (Anderson \cite
{anderson1994},
Andrews \cite{andrews1991} and
Hamilton \cite{hamilton1994} provide the background and numerous
references). In Sections \ref{s:change} and \ref{s:flm-si}, we
illustrate the application of the results of Sections \ref{ss:conv}
and \ref{s:lr} on two problems of recent interest: change point
detection for functional data and the estimation of kernel in the
functional linear model. We show that the detection procedure of
Berkes et al. \cite{berkesgabryshorvathkokoszka2009} must be modified
if the
data exhibit dependence, but the estimation procedure of
Yao, M{\" u}ller and Wang \cite{yaomullerwang2005} is robust to mild
dependence.
Section \ref{s:change} also contains a small simulation study and
a data example. The proofs are collected in the \hyperref[app]{Appendix}.

%s2 ###
\section{Approximable functional
time series} \label{ss:wd}

The notion of weak dependence has, over the past
decades, been formalized in many ways. Perhaps the most
popular are various mixing conditions (see
Doukhan \cite{doukhan1994}, Bradley \cite{bradley2007}), but in recent
years several other approaches have also been introduced
(see Doukhan and Louhichi \cite{doukhanlouhichi1999}
and Wu \cite{wu2005}, \cite{wu2007}, among others).
In time series analysis, moment based measures of dependence, most
notably autocorrelations and cumulants, have gained a universal acceptance.
The measure we consider below is a moment-type quantity, but it is
also related to the mixing conditions as it considers $\sigma$-algebras
$m$ time units apart, with $m$ tending to infinity.

A most direct relaxation of independence is the $m$-dependence.
Suppose $\{X_n \}$ is a sequence of random elements taking
values in a measurable space $S$. Denote by
${\mathcal F}^-_{k} = \sigma\{ \ldots, X_{k-2}, X_{k-1}, X_k\}$ and
${\mathcal F}^+_{k} = \sigma\{X_k, X_{k+1}, X_{k+2}, \ldots\}$
the $\sigma$-algebras generated by the observations up to
time $k$ and after time $k$, respectively. Then
the sequence $\{X_n \}$ is said to be $m$-dependent
if for any $k$, the $\sigma$-algebras ${\mathcal F}^-_{k}$ and
${\mathcal F}^+_{k+m}$
are independent.

Most time series models are not $m$-dependent. Rather, various
measures of dependence decay sufficiently fast, as the distance $m$
between the $\sigma$-algebras ${\mathcal F}^-_{k}$ and ${\mathcal
F}^+_{k+m}$ increases.
However, $m$-dependence can be used as a tool to study properties of
many nonlinear sequences (see, e.g., H{\"o}rmann \cite{hormann2008} and
Berkes, H\"ormann and Schauer \cite{berkeshormannschauer2009} for
recent applications). The
general idea is to approximate $\{X_n, n\in\mathbb{Z}\}$ by $m$-dependent
processes $\{X_{n}^{(m)}, n\in\mathbb{Z}\}$,
$m\geq1$. The goal is to establish that for every $n$
the sequence $\{X_n^{(m)}, m\geq1\}$ converges in some sense
to $X_n$, if we let $m\to\infty$. If the convergence is fast enough,
then one can obtain the limiting behavior
of the original process from
corresponding results for $m$-dependent sequences.
Definition \ref{d:siegi} formalizes this
idea and sets up the necessary framework for the construction of
such $m$-dependent approximation sequences.
The idea of approximating scalar sequences by $m$-dependent
nonlinear moving averages appears already in Section 21 of
Billingsley \cite{billingsley1968}, and it was
developed in several directions by P{\"o}tscher and Prucha \cite
{potscherpruha1997}.
%A version of Definition \ref{d:siegi} for vector
%valued processes was used in \cite{auehormannhorvathreimherr2009}.

In the sequel we let ${H}=L^2([0,1],\mathcal{B}_{[0,1]},\lambda)$
be the Hilbert space of square integrable functions defined on
$[0,1]$. For $f\in{H}$ we set $\|f\|^2=\int_0^1 |f(t)|^2\,dt$.
All our random elements are assumed to
be defined on some common probability space $(\Omega,\mathcal{A},P)$.
For $p\geq1$
we denote by $L^p=L^p(\Omega,\mathcal{A},P)$ the space of
(classes of) real valued random variables such that $\|X\|
_p=(E|X|^p)^{1/p}<\infty$.
Further we let $L_H^p=L_H^p(\Omega,\mathcal{A},P)$ be the space of
$H$ valued random variables $X$ such that
$%\beq\label{e:nu-p}
\nu_p(X) = ( E\|X\|^p )^{1/p}<\infty.
$%\eeq
\begin{definition} \label{d:siegi}
A sequence $\{X_n \}\in L_H^p$ is called
\textit{$L^p$--$m$-approximable}
if each $X_n$ admits the representation,
%
%e2.1 ###
%
\begin{equation}
\label{e:nlma}
X_n = f(\varepsilon_n, \varepsilon_{n-1}, \ldots),
\end{equation}
where the $\varepsilon_i$ are i.i.d. elements taking values in
a measurable space $S$, and $f$ is a measurable function
$f\dvtx S^\infty\to{H}$. Moreover we assume that
if $\{\varepsilon_i'\}$ is an independent
copy of $\{\varepsilon_i\}$ defined on the same probability space,
then letting
%
%e2.2 ###
%
\begin{equation}\label{e:appr}
X_n^{(m)} = f(\varepsilon_n, \varepsilon_{n-1}, \ldots, \varepsilon
_{n-m+1},\varepsilon_{n-m}',
\varepsilon_{n-m-1}',\ldots),
\end{equation}
we have
%
%e2.3 ###
%
\begin{equation}\label{e:mL4}
\sum_{m=1}^\infty\nu_p \bigl(X_m - X_m^{(m)} \bigr) < \infty.
\end{equation}
\end{definition}

For our applications, choosing $p=4$ will be convenient, but any $p\geq1$
can be used, depending on what is needed.
(Our definition makes even sense if $p<1$, but then $\nu_p$ is no longer a norm.)
Definition \ref{d:siegi} implies that $\{X_n\}$ is strictly stationary.
It is clear from the representation
of $X_n$ and $X_n^{(m)}$ that $E\|X_m - X_m^{(m)}\|^p=E\|X_1 -
X_1^{(m)}\|^p$,
so that condition (\ref{e:mL4}) could be formulated solely in terms
of $X_1$ and the approximations $X_1^{(m)}$. Obviously the sequence
$\{X_n^{(m)}, n\in\mathbb{Z}\}$ as defined in (\ref{e:appr})
is \textit{not} $m$-dependent.
To this
end we need to define for each $n$ an independent copy $\{\varepsilon
_k^{(n)}\}
$ of
$\{\varepsilon_k\}$ (this can always be achieved by enlarging the
probability space)
which is then used instead of $\{\varepsilon_k'\}$ to construct $X_n^{(m)}$;
that is,
we set
%
%e2.4 ###
%
\begin{equation}\label{e:apprC}
X_n^{(m)} = f\bigl(\varepsilon_n, \varepsilon_{n-1}, \ldots, \varepsilon
_{n-m+1},\varepsilon{}^{(n)}_{n-m},
\varepsilon_{n-m-1}^{(n)},\ldots\bigr).
\end{equation}
We
will call this method the \textit{coupling construction}.
Since this modification leaves condition (\ref{e:mL4}) unchanged,
we will assume from now on that the $X_n^{(m)}$ are defined by
(\ref{e:apprC}). Then,
for each $m\geq1$, the
sequences $\{X_n^{(m)}, n\in\mathbb{Z}\}$ are strictly stationary
and $m$-dependent, and each $X_n^{(m)}$ is equal in
distribution to $X_n$.

The coupling construction is only
one of a variety of possible $m$-dependent approximations.
In most applications, the measurable space $S$ coincides
with ${H}$, and the $\varepsilon_n$ represent model errors.
In this case, we can set
%
%e2.5 ###
%
\begin{equation}\label{e:trunc}
\tilde X_n^{(m)} = f(\varepsilon_n, \varepsilon_{n-1}, \ldots,
\varepsilon_{n-m+1}, 0, 0,
\ldots).
\end{equation}
The sequence $\{\tilde X_n^{(m)}, n\in\mathbb{Z}\}$ is
strictly stationary and $m$-dependent, but
$X_n^{(m)}$ is no longer equal in
distribution to $X_n$. This is not a big problem but requires
additional lines in the proofs. For the \textit{truncation construction}
(\ref{e:trunc}), condition (\ref{e:mL4}) is replaced by
%
%e2.6 ###
%
\begin{equation}\label{e:mL4T}
\sum_{m=1}^\infty\nu_p \bigl(X_m - \tilde X_m^{(m)} \bigr)< \infty.
\end{equation}
Since $E\|\tilde X_m^{(m)}-X_m^{(m)}\|^p =E\|\tilde X_m^{(m)}-X_m\|^p$,
(\ref{e:mL4T}) implies (\ref{e:mL4}), but not vice versa.
Thus the
coupling construction allows to study a slightly broader
class of time series.

An important question that needs to be addressed at this point
is how our notion of weak dependence compares to other existing
ones. The closest relative of
$L^p$--$m$-approximability is the notion of $L^p$-approximability
studied by P{\"o}tscher and Prucha \cite{potscherpruha1997} for scalar and
vector-valued processes. Since our definition applies with an
obvious modification to sequences with values in any normed vector
spaces $H$ (especially $\mathbb{R}$ or $\mathbb{R}^n$), it can been
seen as a generalization of $L^p$-approximability. There are,
however, important differences. By definition, $L^p$-approximability
only allows for approximations that are, like the truncation
construction, measurable with respect to a finite selection of basis
vectors, $\varepsilon_n,\ldots,\varepsilon_{n-m}$, whereas the coupling
construction does not impose this condition. On the other hand,
$L^p$-approximability is not based on independence of the innovation
process. Instead independence is relaxed to certain mixing
conditions.
Clearly, $m$-dependence implies the CLT, and so our $L^p$--$m$-approximability
implies central limit theorems for practically all
important time series models. As
we have shown in previous papers
\cite{auehormannhorvathreimherr2009,berkeshormannschauer2009,hormann2008,hormann2009}, a scalar version
of this notion has much more
potential than solely giving central limit theorems.
%For functional time series, not many rigorously defined
%models exist yet, but we hope that
%our account can contribute to the development of non-linear functional
%time series by providing the necessary theory.

The concept of weak dependence introduced
in Doukhan and Louhichi \cite{doukhanlouhichi1999}
is defined for scalar
variables in a very general framework
and has been successfully used to prove (empirical)
FCLTs. Like our approach, it does
not require smoothness conditions. Its extensions to problems of
functional data analysis have not been studied yet.

Another approach to weak dependence is a martingale approximation, as
developed in Gordin \cite{gordin1969} and Philipp and Stout
\cite{philippstout1975}. In the context of sequences $\{X_k\}$ of the
form (\ref{e:nlma}), particularly complete results have been proved by
Wu \cite{wu2005,wu2007}. Again, $L^p$--$m$-approximability cannot be
directly compared to approximating martingale conditions; the latter
hold for a very large class of processes, but, unlike
$L^p$--$m$-approximability, they apply only in the context of partial
sums.

The classical approach
to weak dependence, developed in the seminal papers of
Rosenblatt \cite{rosenblatt1956} and Ibragimov \cite{ibragimov1962},
uses the strong mixing property and its variants like
$\beta$, $\phi$, $\rho$ and $\psi$ mixing. The general
idea is to measure the maximal dependence between two events
lying in the ``past'' $\mathcal{F}_k^-$
and in the ``future'' $\mathcal{F}_{k+m}^+$, respectively.
The fading memory is described by this maximal dependence
decaying to
zero for $m$ growing to $\infty$.
For example, the $\alpha$-mixing
coefficient is given by
\[
\alpha_m=\sup\{|P(A\cap B)-P(A)P(B)|A\in\mathcal{F}_k^-, B\in
\mathcal{F}_{k+m}^+, k\in\mathbb{Z}\}.
\]
A sequence is called $\alpha$-mixing (strong mixing) if $\alpha_m\to0$
for $m\to\infty$.

This method yields very sharp
results (for a complete account of the classical theory (see Bradley
\cite{bradley2007}), but
verifying mixing conditions of the above type is not easy,
whereas the verification of $L^p$--$m$-approximability
is almost immediate as our examples below show.
This is because the $L^p$--$m$-approximability condition uses
directly the model specification $X_n = f(\varepsilon_n, \varepsilon
_{n-1}, \ldots)$.
Another problem is that even when mixing
applies (e.g., for Markov processes), it typically requires
strong smoothness conditions. For example, for the AR(1) process
\[
Y_k=\tfrac{1}{2} Y_{k-1}+\varepsilon_k
\]
with Bernoulli innovations, strong mixing fails to hold (cf.
Andrews \cite{andrews1984}). Since
$c$-mixing, where $c$ is either of $\psi$, $\phi$, $\beta$ or $\rho$,
implies $\alpha$-mixing, $\{Y_k\}$ above satisfies none of these
mixing conditions, whereas Example \ref{ex:m-far}
shows that the AR(1) process is $L^p$--$m$-approximable
without requiring any smoothness properties for the
innovations process.
Consequently our condition does not imply strong mixing.
On the other hand, $L^p$--$m$-approximability
is restricted to a more limited class of processes, namely processes
allowing the representation $X_n = f(\varepsilon_n, \varepsilon
_{n-1}, \ldots)$.
We emphasize, however, that all time series
models used in practice (scalar, vector or functional)
have this representation
(cf. \cite{priestley1988,stine1997,tong1990}), as an immediate consequence of
their ``forward'' dynamics, for example, their
definitions by a stochastic recurrence equations.
See the papers of Rosenblatt
\cite{rosenblatt1959,rosenblatt1961,rosenblatt1971}
for sufficient criteria.

We conclude that \textit{$L^p$--$m$-approximability
is not directly comparable with
classical mixing coefficients.}

The following lemma shows how $L^p$--$m$-approximability is
unaffected by linear transformations, whereas independence assumptions
are needed for product type operations.
\begin{lemma}\label{l:mani}
Let $\{X_n\}$ and $\{Y_n\}$ be two $L^p$--$m$-approximability
sequences in $L_H^p$. Define:
\begin{itemize}
\item$Z_n^{(1)}=A(X_n)$, where $A\in\mathcal{L}$;
\item$Z_n^{(2)}=X_n+Y_n$;
\item$Z_n^{(3)}=X_n\circ Y_n$ (point-wise multiplication);
\item$Z_n^{(4)}=\langle X_n,Y_n \rangle$;
\item$Z_n^{(5)}=X_n\otimes Y_n$.
\end{itemize}
Then $\{Z_n^{(1)}\}$ and $\{Z_n^{(2)}\}$ are
$L^p$--$m$-approximable
sequences in $L_H^p$.
If $X_n$ and $Y_n$ are independent, then
$\{Z_n^{(4)}\}$ and $\{Z_n^{(5)}\}$ are
$L^p$--$m$-approximable sequences in the respective
spaces.
If\vspace*{1pt}
$E\sup_{t\in[0,1]}|X_n(t)|^p+E\sup_{t\in[0,1]}|Y_n(t)|^p<\infty$,
then
$\{Z_n^{(3)}\}$ is $L^p$--$m$-approximable in $L_H^p$.
\end{lemma}
\begin{pf}
The first two relations are immediate. We exemplify the
rest of the
simple proofs for $Z_n=Z_n^{(5)}$. For this we set
$Z_m^{(m)}=X_m^{(m)}\otimes Y_m^{(m)}$ and note that $Z_m$ and $Z_m^{(m)}$
are (random) kernel operators, and thus Hilbert--Schmidt operators.
Since
\begin{eqnarray*}
\bigl\|Z_m-Z_m^{(m)}\bigr\|_\mathcal{L}&\leq&\bigl\|Z_m-Z_m^{(m)}\bigr\|_\mathcal{S}\\
&\leq&\biggl(\dint\bigl(X_m(t)Y_m(s)-X_m^{(m)}(t)Y_m^{(m)}(s)
\bigr)^2\,dt\,ds \biggr)^{1/2}\\
&\leq&\sqrt{2} \bigl(\|X_m\|\bigl\|Y_m-Y_m^{(m)}\bigr\|+\bigl\|Y_m^{(m)}\bigr\|\bigl\|
X_m-X_m^{(m)}\bigr\| \bigr),
\end{eqnarray*}
the proof follows from the independence of $X_n$ and $Y_n$.
\end{pf}

The proof shows that our assumption can be modified and independence is
not required. However, if $X,Y$ are not independent, then $E|XY|\neq E|X|E|Y|$.
We have then to use the Cauchy--Schwarz inequality
and obviously need $2p$ moments.

We want to point out that only a straightforward modification
is necessary in order to generalize the theory of this paper to
noncausal processes $ X_n=f(\ldots,\varepsilon_{n+1},\varepsilon
_n,\varepsilon
_{n-1},\ldots).
$ Our framework
can be also extended to nonstationary sequences, for example, those
of the
form (\ref{e:nlma}) where
$\{\varepsilon_k\}$ is a sequence of independent, but not necessarily
identically
distributed, or random variables where
\[
X_n=f_n(\varepsilon_n,\varepsilon_{n-1},\ldots).
\]
The $m$-dependent coupled process can be defined in
the exact same way as in the stationary case
\[
X_n^{(m)}=f_n\bigl(\varepsilon_n,\varepsilon_{n-1},\ldots,\varepsilon_{n-m+1},
\varepsilon_{n-m}^{(n)},\varepsilon_{n-m-1}^{(n)},\ldots\bigr).
\]
A generalization of our method to nonstationarity would be
useful, especially when the goal is to develop
methodology for locally stationary data.
Such work is, however, beyond the intended scope of this paper.

We now illustrate the applicability of Definition \ref{d:siegi} with
several examples. Let
$\mathcal{\mathcal{L}}=\mathcal{\mathcal{L}}({H},{H})$ be the set of
bounded linear operators from ${H}$ to ${H}$. For $A\in
\mathcal{\mathcal{L}}$ we define the operator norm
$\|A\|_{\mathcal{L}}=\sup_{\|x\|\leq1}\|Ax\|$.
If the operator is
Hilbert--Schmidt, then we denote with $\|A\|_\mathcal{S}$ its
Hilbert--Schmidt norm. Recall that for any Hilbert--Schmidt
operator $A\in\mathcal{L}$, $\|A\|_\mathcal{L}\leq\|A\|_{\mathcal S}$.
\begin{example}[(Functional autoregressive process)]\label{ex:m-far}
Suppose $\Psi\in\mathcal{L}$ satisfies
$\|\Psi\|_{\mathcal{L}}<1$.
Let
$\varepsilon_n\in L_H^2$ be i.i.d. with mean zero.
Then
there is a unique stationary sequence of random elements
$X_n \in L_H^2$ such that
%
%e2.7 ###
%
\begin{equation}\label{e:ar1-si}
X_n(t) = \Psi(X_{n-1})(t)+\varepsilon_n(t).
\end{equation}
For details see Chapter 3 of Bosq \cite{bosq2000}.
%Specifically if $\Psi$ is an integral Hilbert--Schmidt operator with
%kernel
%$\psi(t, s), \ t,s \in[0,1]$, then
%$$
%X_n(t) = \int\psi(t,s) X_{n-1}(s)\,ds + \eg_n(t),
%$$
%and since
%$$
% \|\Psi\|_{\mathcal{L}}^2 \leq\|\Psi\|_{\cS}^2= \dint\psi^2(t,
%s)dt\,ds,
%$$
%a sufficient condition for the existence
%of a unique stationary and causal solution is
%$\dint\psi^2(t, s)dt\,ds<1$.
The AR(1) sequence (\ref{e:ar1-si}) admits the expansion
$
X_n = \sum_{j=0}^\infty\Psi^j(\varepsilon_{n-j})
$
where $\Psi^j$ is the $j$th iterate of the operator $\Psi$.
We thus set
$
X_n^{(m)} = \sum_{j=0}^{m-1} \Psi^j(\varepsilon_{n-j})
+\sum_{j=m}^{\infty}\Psi^j(\varepsilon_{n-j}^{(n)}).
$
It is easy to verify that for every $A$ in $\mathcal{L}$,
$
\nu_p(A(Y)) \le\|A\|_{\mathcal{L}} \nu_p(Y).
$
Since
$
X_m - X_m^{(m)}
= \sum_{j=m}^\infty(\Psi^j(\varepsilon_{m-j})-\Psi
^j(\varepsilon
_{m-j}^{(m)}) ),
$
it follows that
$
\nu_p(X_m - X_m^{(m)})
\le\break {2\sum_{j=m}^\infty}\|\Psi\|_\mathcal{L}^j \nu_p(\varepsilon_0)
= O(1) \times\nu_p(\varepsilon_0) \|\Psi\|_{\mathcal{L}}^m.
$
By assumption $\nu_2(\varepsilon_0)<\infty$ and therefore
$
\sum_{m=1}^\infty\nu_2(X_m - X_m^{(m)}) < \infty,
$
so condition (\ref{e:mL4T}) holds with $p\geq2$, as long
as $\nu_p(\varepsilon_0)<\infty$.
\end{example}

The argument in the above example shows
that a sufficient condition to obtain $L^p$--$m$-approximability
is
\begin{eqnarray*}
&&\|f(a_m,\ldots,a_1,x_0,x_{-1},\ldots)-f(a_m,\ldots
,a_1,y_0,y_{-1},\ldots)\|\\
&&\qquad \leq c_m \|f(x_0,x_{-1},\ldots)-f(y_0,y_{-1},\ldots)\|,
\end{eqnarray*}
where $\sum_{m\geq1} c_m<\infty$. This holds for a functional
AR(1) process and offers an attractive sufficient and
distribution-free condition for $L^p$--$m$-approxim\-ability.
The interesting
question, whether one can impose some other, more general conditions
on the function $f$ that would imply $L^p$--$m$-approxim\-ability
remains open.
For example, the simple criterion above does not apply to general linear
processes.
We recall that a sequence $\{ X_n \}$ is said to be a \textit{linear process
in ${H}$} if
$
X_n = \sum_{j=0}^\infty\Psi_j(\varepsilon_{n-j})
$
where the errors $\varepsilon_n\in L^2_H$ are i.i.d. and zero mean,
and each
$\Psi_j$
is a bounded operator. If
$
\sum_{j=1}^\infty\|\Psi_j\|_\mathcal{L}^2 < \infty,
$
then the series defining $X_n$ converges a.s. and in $L_H^2$
(see Section 7.1 of Bosq \cite{bosq2000}).

A direct verification,
following the lines of Example \ref{ex:m-far}, yields sufficient
conditions for a general linear
process to be $L^p$--$m$-approximable.
\begin{proposition} \label{p:m-lin}
Suppose $\{ X_n \}\in L_H^2$ is a linear process whose errors satisfy
$\nu_p(\varepsilon_0)< \infty$, $p\geq2$. The operator coefficients satisfy
%
%e2.8 ###
%
\begin{equation}\label{e:Lpafg}
{\sum_{m=1}^\infty\sum_{j=m}^\infty}\|\Psi_j\|< \infty.
\end{equation}
Then $\{ X_n \}$ is $L^p$--$m$-approximable.
\end{proposition}

We note that condition (\ref{e:Lpafg}) is comparable
to the usual assumptions made in the scalar case.
For a scalar linear process the weakest
possible condition for weak dependence is
%
%e2.9 ###
%
\begin{equation}\label{e:sm}
{\sum_{j=0}^\infty}|\psi_j|<\infty.
\end{equation}
If it is violated, the
resulting time series are referred to as strongly dependent, long memory,
long-range dependent or persistent.
Recall that (\ref{e:sm}) merely ensures the existence
of fundamental population objects like an absolutely summable
autocovariance sequence or a bounded spectral density. It is, however,
too weak to establish any statistical results.
For example, for the asymptotic normality of the sample
autocorrelations we need $\sum j \psi_j^2 < \infty$, for
the convergence of the periodogram ordinates
$\sum\sqrt{j} |\psi_j|< \infty$.
Many authors assume $ \sum{j} |\psi_j|< \infty$ to be able to use
all these basic results.
The condition $\sum{j} |\psi_j|< \infty$
is equivalent to (\ref{e:Lpafg}).

We next give a simple example of a nonlinear $L^p$--$m$-approximable
sequence. It is based on the model used by
Maslova et al. \cite{maslovakokoszkasz2009b} to simulate the so-called solar
quiet (Sq) variation in magnetometer records (see also
Maslova et al. \cite{maslovakokoszkasz2009a}). In that model, $X_n(t) =
U_n(S(t) + Z_n(t))$ represents the part of the magnetometer record on
day $n$ which reflects the magnetic field generated by ionospheric
winds of charged particles driven by solar heating. These winds flow
in two elliptic cells, one on each day-side of the equator. Their
position changes from day to day, causing
a different appearance of the curves, $X_n(t)$, with changes in the
amplitude being most pronounced. To simulate this behavior, $S(t)$ is
introduced as the typical pattern for a specific magnetic observatory, $Z_n(t)$,
as the change in shape on day $n$ and the scalar random variable
$U_n$ as the amplitude on day $n$. With this motivation, we formulate
the following example.
\begin{example}[(Product model)] \label{ex:sq}
Suppose $\{ Y_n \}\in L_H^p$ and $\{ U_n \}\in L^p$ are both
$L^p$--$m$-approximable sequences, independent of each other.
The respective representations are
$Y_n = g(\eta_1, \eta_2, \ldots)$ and
$U_n = h(\gamma_1, \gamma_2, \ldots)$.
Each of these sequences could be a linear sequence
satisfying the assumptions of Proposition \ref{p:m-lin}, but
they need not be. The sequence $X_n(t)= U_n Y_n(t)$ is then a nonlinear
$L^p$--$m$-approximable sequence with the underlying
i.i.d. variables $\varepsilon_n = (\eta_n, \gamma_n)$. This follows by
after a slight modification from Lemma \ref{l:mani}.
%To see this, set
%$X_n^{(m)}(t) = U_n^{(m)} Y_n^{(m)}(t)$ and observe that
%$
%$
%Using the independence of $\{ Y_n \}$ and $\{ U_n \}$ it can be easily
%shown that
%%\[
%%\nu_p((U_m - U_m^{(m)})Y_m)
%%= \{ E [\int(U_m - U_m^{(m)})^2 Y_m^2(t) dt ]^{p/2} \}^{1/p}
%%\]
%%\[
%%= \{ E|U_m - U_m^{(m)}|^p E[\int Y_m^2(t) dt ]^{p/2}\}^{1/p}
%%= \|U_m - U_m^{(m)}\|_p \nu_p(Y_0)
%%\]
%$
%= \|U_m - U_m^{(m)}\|_p \nu_p(Y_0)
%$
%and
%$
%= \|U_0\|_p \nu_p( Y_m - Y_m^{(m)} ).
%$
\end{example}

Example \ref{ex:sq} illustrates the principle that in order for
products of $L^p$--$m$-approximable sequences to be
$L^p$--$m$-approximable, independence must be assumed. It does not
have to be assumed as directly as in Example \ref{ex:sq}; the
important point being that appropriately-defined functional Volterra
expansions should not contain diagonal terms so that moments do not
pile up. Such expansions exist (see, e.g.,
Giraitis, Kokoszka and Leipus \cite{giraitiskokoszkaleipus2000ET}, for
all nonlinear scalar
processes used to model financial data). The model $X_n(t)= Y_n(t) U_n$
is similar to the popular scalar stochastic volatility model $r_n =
v_n \varepsilon_n$ used to model returns $r_n$ on a speculative asset. The
dependent sequence $\{v_n\}$ models volatility, and the i.i.d. errors
$\varepsilon_n$, independent of the $v_n$, generate unpredictability in
returns.

Our next examples focus on functional extensions of popular
nonlinear models, namely the bilinear model of
\cite{grangeranderson1978}
and the celebrated ARCH model of Engle~\cite{engle1982}.
Both models will be treated in more detail in
forthcoming papers. Proofs of Propositions \ref{p:bilinear}
and \ref{p:farch} are available upon request.
\begin{example}[(Functional bilinear process)]\label{ex:bilinear}
Let $(\varepsilon_n)$ be an $H$-valued i.i.d. sequence and let
$\psi\in H\otimes H$ and $\phi\in H\otimes H\otimes H$. Then the
process defined as the recurrence equation,
\[
X_{n+1}(t) = \int\psi(t,s) X_n(s) \,ds
+\dint\phi(t,s,u) X_n(s)\varepsilon_n(u)\,ds\,du + \varepsilon_{n+1}(t),
\]
is called \textit{functional bilinear process}.

A neater notation
can be achieved by defining
$\psi\dvtx H\to H$, the kernel operator with the kernel function
$\phi(t,s)$, and $\phi_n\dvtx H\to H$, the random kernel operator
with kernel
\[
\phi_n(t,s)=\int\phi(t,s,u) \varepsilon_n(u)\,du.
\]
In this notation, we have
%
%e2.10 ###
%
\begin{equation}\label{e:bilinear}
X_{n+1}=(\psi+\phi_n)(X_n)+\varepsilon_{n+1}
\end{equation}
with the usual convention that $(A+B)(x)=A(x)+B(x)$
for operators $A,B$.
The product of two operators
$AB(x)$ is interpreted as successive application $A(B(x))$.

%Then iteration of \eqref{e:bilinear} gives
%$$
%X_{n+1}=\sum_{k=0}^m\prod_{j=0}^{k-1}(\psi+\phi_{n-j})(
%$$
%where we define $\prod_{j=0}^{-1}(\psi+\phi_{n-j})(x)=x$.

A formal solution to (\ref{e:bilinear}) is
%
%e2.11 ###
%
\begin{equation}\label{e:bilinearsolution}
X_{n+1}=\sum_{k=0}^\infty\prod_{j=0}^{k-1}
(\psi+\phi_{n-j})(\varepsilon_{n+1-k})
\end{equation}
and the approximating sequence is defined by
\[
\tilde{X}_n^{(m)}=
\sum_{k=0}^m\prod_{j=0}^{k-1}(\psi+\phi_{n-j})(\varepsilon_{n+1-k}).
\]

The following proposition establishes sufficient conditions for
the $L^p$--$m$-approximability.
\end{example}
\begin{proposition}\label{p:bilinear}
Let $\{ X_n \}$ be
the functional bilinear process defined in (\ref{e:bilinear}).
If $E\log\|\varepsilon_0\|<\infty$ and $E\log\|\psi+\phi_0\|<0$, then
a unique strictly stationary solution for this equation exists.
The solution has ($L^2$-)representation (\ref{e:bilinearsolution}). If
$\nu_p ((\psi+\phi_{0})(\varepsilon_{0}) )<\infty$ and
$E\|\psi+\phi_{0}\|_\mathcal{S}^p<1$, the process is
$L^p$--$m$-approximable.
\end{proposition}
\begin{example}[(Functional ARCH)]\label{ex:ARCH}
Let $\delta\in H$ be a positive function and
let $\{\varepsilon_k\}$ an i.i.d.
sequence in $L_H^4$.
Further, let $\beta(s,t)$ be a nonnegative kernel function in $L^2([0,1]^2,
\mathcal{B}_{[0,1]}^2,\lambda^2)$.
Then we call the process
%
%e2.12 ###
%
\begin{equation}\label{e:y}
y_k(t)=\varepsilon_k(t)\sigma_k(t),\qquad t\in[0,1],
\end{equation}
where
%
%e2.13 ###
%
\begin{equation}\label{e:sigma}
\sigma_k^2(t)=\delta(t)+\int_0^1\beta(t,s) y_{k-1}^2(s)\,ds,
\end{equation}
the \textit{functional $\mathrm{ARCH}(1)$ process}.

Proposition \ref{p:farch} establishes conditions for the existence of
a strictly stationary solution to (\ref{e:y}) and
(\ref{e:sigma}) and its $L^p$--$m$-approximability.
\end{example}
\begin{proposition} \label{p:farch}
Define $K(\varepsilon_1^2)= (\dint\beta^2(t,s)\varepsilon
_1^4(s) \,ds\,dt )^{1/2}$.
If there is some $p>0$ such that
$E \{K(\varepsilon_1^2) \}^p<1$ then (\ref{e:y})
and (\ref{e:sigma})
have a unique strictly stationary and causal solution
and the sequence $\{ y_k \}$ is\break $L^p$--$m$-approximable.
\end{proposition}

%s3 ###
\section{Convergence of eigenvalues
and eigenfunctions}\label{ss:conv}

Denote by $C=E[\langle X$,\break $\cdot\rangle X]$ the covariance operator of some
$X\in
L_H^2$. The eigenvalues and eigenfunctions of $C$ are a fundamental
ingredient for principal component analysis which is a key technique
in functional data analysis. In practice, $C$ and its
eigenvalues/eigenfunctions are unknown and must be estimated. The
purpose of this section is to prove consistency of the corresponding
estimates for $L^4$-$m$-approximable sequences. The results derived
below will be applied in the following sections. We start with some
preliminary results.

Consider two compact operators $C, K \in{\mathcal L}$ with singular
value decompositions
%
%e3.1 ###
%
\begin{equation}\label{e:C-K}
C(x) = \sum_{j=1}^\infty\lambda_j \langle x, v_j \rangle f_j,\qquad
K(x) = \sum_{j=1}^\infty\gamma_j \langle x, u_j \rangle g_j.
\end{equation}
The following lemma is proven in Section VI.1 of
(see Gohberg, Golberg and Kaashoek \cite{gohberg1990}, Corollary 1.6,
page 99).
\begin{lemma} \label{l:B4.2}
Suppose $C, K \in{\mathcal L}$ are two compact operators with singular
value decompositions (\ref{e:C-K}). Then, for each $j\ge1$,
$
|\gamma_j - \lambda_j| \le\|K-C\|_{\mathcal L}.
$
\end{lemma}

We now tighten the conditions on the operator $C$ by assuming
that it is Hilbert--Schmidt, symmetric and positive definite.
These conditions imply that
$f_j=v_j$ in~(\ref{e:C-K}), $C(v_j)=\lambda_jv_j$
and $\sum_j \lambda_j^2<\infty$. Consequently $\lambda_j$ are eigenvalues
of $C$ and $v_j$ the corresponding eigenfunctions. We also define
\[
v_j^\prime={\hat c}_j v_j,\qquad {\hat c}_j=\operatorname{sign}(\langle
u_j,v_j \rangle).
\]
Using Lemma \ref{l:B4.2}, the next lemma can be established by
following the lines of the proof of Lemma 4.3 of Bosq \cite{bosq2000}.
\begin{lemma} \label{l:B4.3}
Suppose $C, K \in{\mathcal L}$ are two compact operators with singular
value decompositions (\ref{e:C-K}).
If $C$ is Hilbert--Schmidt, symmetric and positive definite,
and its eigenvalues satisfy
%
%e3.2 ###
%
\begin{equation}\label{e:la>d}
\lambda_1 > \lambda_2> \cdots> \lambda_d > \lambda_{d+1},
\end{equation}
then
\[
\|u_j - v_j^\prime\| \le\frac{2\sqrt{2}}{\alpha_j} \|K-C\|
_{\mathcal L},\qquad
1 \le j\le d,
\]
where $\alpha_1 = \lambda_1- \lambda_2$ and
$
\alpha_j = \min( \lambda_{j-1} - \lambda_j, \lambda_j - \lambda
_{j+1}), 2 \le j\le d.
$
\end{lemma}

Let $\{ X_n \}\in L_H^2$ be a
stationary sequence with covariance
operator $C$. In principle we could now develop a general theory
for $H$ valued sequences, where $H$ is an arbitrary separable Hilbert space.
In practice, however, the case $H=L^2([0,1],\mathcal
{B}_{[0,1]},\lambda)$
is most important.
In order to be able to fully use the structure of $H$ and
and not to deal with technical assumptions,
we need the two basic regularity conditions
below, which will be assumed throughout the paper without further notice.
\begin{assumption}\label{a:reg}
(i) Each $X_n$ is measurable
$(\mathcal{B}_{[0,1]}\times\mathcal{A})/\mathcal{B}_\mathbb
{R}$.

(ii) $\sup_{t\in[0,1]}E|X(t)|^2<\infty$.
\end{assumption}

Assumption \ref{a:reg}(i) is necessary in order that the sample
paths of $X_n$ are measurable. Together with (ii) it also implies
that $C$ is an integral operator with kernel
$
c(t,s) =\operatorname{Cov}(X_1(t),X_1(s))
$
whose estimator is
%
%e3.3 ###
%
\begin{equation}\label{e:estkern}
{\hat c}(t,s) = N^{-1} \sum_{n=1}^N \bigl(X_n(t)-\bar{X}_N(t)\bigr)
\bigl(X_n(s)-\bar{X}_N(s)\bigr).
\end{equation}
Then natural estimators of the eigenvalues $\lambda_j$
and eigenfunctions $v_j$ of $C$ are the
eigenvalues ${\hat\lambda}_j$ and eigenfunctions ${\hat v}_j$ of
$\hat{C}$,
the operator with the kernel
(\ref{e:estkern}). By Lemmas \ref{l:B4.2} and \ref{l:B4.3} we can
bound the estimation errors for eigenvalues and eigenfunctions by
$\|C-\hat{C}\|_{{\mathcal S}}^2$.
Mas and Mennetau \cite{masmennetau2003}
show that transferring asymptotic results
from the operators to the eigenelements holds quite generally,
including a.s. convergence, weak convergence or
large deviation principles.
This motivates the next result.
\begin{theorem} \label{t:ChatW}
Suppose $\{ X_n \}\in L_H^4$ is an $L^4$--$m$-approximable
sequence with covariance operator $C$. Then
there is some constant $U_X<\infty$, which does not depend on $N$,
such that
%
%e3.4 ###
%
\begin{equation}\label{e:s2}
NE\|\hat C - C\|_{{\mathcal S}}^2 \le U_X.
\end{equation}
If the $X_n$ have zero mean, then we can choose
%
%e3.5 ###
%
\begin{equation}\label{e:UX}
U_X= \nu_4^4(X) + 4\sqrt{2} \nu_4^3(X)
\sum_{r=1}^\infty\nu_4\bigl(X_r - X_r^{(r)}\bigr).
\end{equation}
\end{theorem}

The proof of Theorem \ref{t:ChatW} is given in Section \ref{s:p-si}.
Let us note that by Lem\-ma~\ref{l:B4.2} and Theorem \ref{t:ChatW},
\[
N E [|\lambda_j -{\hat\lambda}_j|^2]\le NE\|\hat C - C\|_{{\mathcal L}}^2
\le NE\|\hat C - C\|_{{\mathcal S}}^2\le U_X.
\]
Assuming (\ref{e:la>d}), by Lemma
\ref{l:B4.3} and Theorem \ref{t:ChatW},
[$\hat{c}_j=\operatorname{sign}(\langle{\hat v}_j,v_j \rangle)$],
\[
N E [\|\hat c_j{\hat v}_j -v_j\|^2]\le\biggl(\frac{2\sqrt{2}}{\alpha
_j}\biggr)^2
NE\|\hat C - C\|_{{\mathcal L}}^2
\le\frac{8}{\alpha_j^2} NE\|\hat C - C\|_{{\mathcal S}}^2
\le\frac{8U_X}{\alpha_j^2}
\]
with the $\alpha_j$ defined in Lemma \ref{l:B4.3}.

These inequalities establish the following result.
\begin{theorem} \label{t:B4W}
Suppose $\{ X_n \}\in L_H^4$ is an $L^4$--$m$-approximable
sequence and assumption (\ref{e:la>d})
holds.
Then, for $1\le j \le d$,
%
%e3.6 ###
%
\begin{equation}\label{e:eigen-si}\qquad
\limsup_{N\to\infty} N E [|\lambda_j -{\hat\lambda}_j|^2]<
\infty,\qquad
\limsup_{N\to\infty} N E [\|\hat c_j{\hat v}_j -v_j\|^2]< \infty.
\end{equation}
\end{theorem}

Relations (\ref{e:eigen-si}) have become a fundamental tool for establishing
asymptotic properties of procedures for functional simple random samples
which are based on the functional principal components. Theorem
\ref{t:B4W} shows that in many cases one can expect that these
properties will remain the same under weak dependence; an important
example is discussed in Section \ref{s:flm-si}.
The empirical covariance kernel (\ref{e:estkern}) is, however, clearly designed
for simple random samples, and may not be optimal for representing
dependent data in the most ``useful'' way. The term ``useful'' depends on
the application. Kargin and Onatski \cite{karginonatski2008} show that
a basis different
than the eigenfunctions $v_k$ is optimal for prediction with a functional
AR(1) model. An interesting open problem is how to construct a
basis optimal in some general sense for dependent data.
In Section \ref{s:lr} we focus on a related, but different, problem
of constructing a matrix which ``soaks up'' the dependence in a manner
that allows the extension of many multivariate time series procedures
to a
functional setting. The construction of this matrix involves
\textit{arbitrary} basis vectors $v_k$ estimated
by ${\hat v}_k$ in such a way that (\ref{e:eigen-si}) holds.

%s4 ###
\section{Estimation of the long-run
variance} \label{s:lr}

The main results of this section are Corollary
\ref{c:lr-f} and Proposition \ref{p:diffbartlett}
which state that the long-run variance matrix obtained by
projecting the data on the functional principal components
can be consistently estimated. The concept of the long-run variance,
while fundamental in time series analysis, has not been studied
for functional data, and not even for scalar approximable sequences.
It is therefore necessary to start with some preliminaries which lead
to our main results and illustrate the role of the
$L^p$--$m$-approximability.

Let $\{X_n\}$ be a scalar (weakly) stationary sequence.
Its long-run variance is defined as
$
\sigma^2=\sum_{j\in\mathbb{Z}}\gamma_j,
$
where $\gamma_j = \operatorname{Cov}(X_0, X_j)$, provided this series
is absolutely convergent.
Our first lemma shows that this is the case for $L^2$--$m$-approximable
sequences.
\begin{lemma} \label{l:cov-sum}
Suppose $\{ X_n \}$ is a scalar $L^2$--$m$-approximable
sequence. Then its autocovariance function $\gamma_j = \operatorname
{Cov}(X_0, X_j)$
is absolutely summable, that is, $\sum_{j=-\infty}^\infty|\gamma_j| <
\infty$.
\end{lemma}
\begin{pf} Observe that for $j>0$,
\[
\operatorname{Cov}(X_0, X_j) =\operatorname{Cov}\bigl(X_0, X_j -
X_j^{(j)}\bigr) + \operatorname{Cov}\bigl(X_0, X_j^{(j)}\bigr).
\]
Since
\[
X_0 = f(\varepsilon_0, \varepsilon_{-1}, \ldots),\qquad
X_j^{(j)} = f^{(j)}\bigl(\varepsilon_j, \varepsilon_{j-1}, \ldots,
\varepsilon_1,\varepsilon
_0^{(j)},\varepsilon_{-1}^{(j)},\ldots\bigr),
\]
the random variables $X_0$ and $X_j^{(j)}$ are independent, so
$\operatorname{Cov}(X_0, X_j^{(j)}) = 0$, and
\[
|\gamma_j| \le[EX_0^2]^{1/2} \bigl[ E \bigl(X_j - X_j^{(j)}\bigr)^2 \bigr]^{1/2}.
\]
\upqed\end{pf}

The summability of the autocovariances is the fundamental property
of weak dependence because then $N \operatorname{Var}[\bar X_N] \to
\sum_{j=-\infty}^\infty\gamma_j$; that is, the variance of the
sample mean converges to zero at the rate $N^{-1}$, the same
as for i.i.d. observations.
A popular approach to the estimation of
the long-run variance is to use the kernel estimator
\[
\hat\sigma^2 = \sum_{|j| \le q} \omega_q(j) \hat\gamma_j,\qquad
\hat\gamma_j
= \frac{1}{N}\sum_{i=1}^{N-|j|} (X_i - \bar X_N) \bigl(X_{i+|j|} - \bar X_N\bigr).
\]
Various weights $\omega_q(j)$ have been
proposed and their optimality properties studied (see
Anderson \cite{anderson1994} and Andrews \cite{andrews1991},
among others). In theoretical work, it is typically assumed that
the bandwith $q$ is a deterministic function of the sample size
such that $q= q(N)\to\infty$ and $q= o(N^r)$, for some $0<r\le1$.
We will use the following assumption:
\begin{assumption} \label{a:wq}
The bandwidth $q= q(N)$ satisfies
$
q \to\infty, q^2/N \to0,
$
and the weights satisfy $\omega_q(j)= \omega_q(-j)$ and
%
%e4.1 ###
%
\begin{equation}\label{e:og-b}
|\omega_q(j)| \leq b
\end{equation}
and, for every fixed $j$,
%
%e4.2 ###
%
\begin{equation}\label{e:og-c}
\omega_q(j) \to1.
\end{equation}
\end{assumption}

All kernels used in practice have symmetric weights
and satisfy conditions (\ref{e:og-b}) and (\ref{e:og-c}).

The absolute summability of the autocovariances is
not enough to establish the consistency of the kernel estimator
$\hat\sigma^2$. Traditionally, summability of the cumulants has been
assumed to control the fourth order structure of the data.
Denoting $\mu= EX_0$, the fourth order cumulant of a stationary
sequence is defined by
\[
\kappa(h,r,s) =
\operatorname{Cov}\bigl((X_0 - \mu)(X_h - \mu), (X_r-\mu)(X_s-\mu)\bigr)
-\gamma_r\gamma_{h-s} - \gamma_s\gamma_{h-r}.
%=E[(X_0- \mu)(X_{h} - \mu) (X_{r} - \mu)(X_{s} - \mu)]
%-(\ga_h\ga_{r-s} + \ga_r\ga_{h-s} + \ga_r\ga_{h-r}).
\]
The ususal sufficient condition for the consistency of $\hat\sigma$ is
%
%e4.3 ###
%
\begin{equation}\label{e:cum-cond}
{\sum_{h=-\infty}^\infty\sum_{r=-\infty}^\infty\sum_{s=-\infty
}^\infty}
|\kappa(h,r,s)|< \infty.
\end{equation}
%
%If $X_n = \sum_{j=0}^\infty c_j \eg_{n-j}$ is a linear sequence
%with i.i.d. finite fourth
%order errors, then
%c_k c_{k+h} c_{k+r} c_{k+s},
%so condition \refeq{cum-cond}
%holds, provided $\sum_j |c_j| < \infty$.
%}%comment
Recently, Giraitis et al. \cite{giraitiskokoszkaleipus2003RV} showed that
condition (\ref{e:cum-cond}) can be replaced by a weaker condition,
%
%e4.4 ###
%
\begin{equation}\label{e:gklt-cond}
{\sup_h \sum_{r=-\infty}^\infty\sum_{s=-\infty}^\infty}
|\kappa(h,r,s)|< \infty.
\end{equation}
A technical condition we need is
%
%e4.5 ###
%
\begin{equation}\label{e:lr-ct}
N^{-1} \sum_{k,l=0}^{q(N)} \sum_{r=1}^{N-1}
\bigl|\operatorname{Cov}\bigl( X_0 \bigl(X_{k}- X_{k}^{(k)}\bigr),
X_r^{(r)} X_{r+\ell}^{(r+\ell)} \bigr)\bigr|\to0.
\end{equation}
By analogy to condition (\ref{e:gklt-cond}), it can be replaced
by a much stronger, but a more transparent condition,
%
%e4.6 ###
%
\begin{equation}\label{e:lr-cs}
\sup_{k,l\ge0} \sum_{r=1}^{\infty}
\bigl|\operatorname{Cov}\bigl( X_0 \bigl(X_{k}- X_{k}^{(k)}\bigr),
X_r^{(r)} X_{r+\ell}^{(r+\ell)} \bigr)\bigr|< \infty.
\end{equation}

To explain the intuition behind conditions (\ref{e:lr-ct})
and (\ref{e:lr-cs}), consider the linear process $X_k = \sum
_{j=0}^\infty
c_j X_{k-j}$.
For $k \ge0$,
\[
X_k - X_{k}^{(k)}
= \sum_{j=k}^\infty c_j \varepsilon_{k-j} -\sum_{j=k}^\infty c_j
\varepsilon_{k-j}^{(k)}.
\]
Thus $X_0 (X_{k}- X_{k}^{(k)})$ depends on
%
%e4.7 ###
%
\begin{equation}\label{e:eg0}
\varepsilon_0, \varepsilon_{-1}, \varepsilon_{-2}, \ldots
\quad\mbox{and}\quad
\varepsilon_0^{(k)}, \varepsilon_{-1}^{(k)}, \varepsilon_{-2}^{(k)},
\ldots
\end{equation}
and $X_r^{(r)} X_{r+\ell}^{(r+\ell)}$ depends on
\[
\varepsilon_{r+\ell}, \ldots, \varepsilon_1,
\varepsilon_0^{(r)} \varepsilon_{-1}^{(r)}, \varepsilon_{-2}^{(r)},
\ldots\quad\mbox
{and}\quad
\varepsilon_0^{(r+\ell)} \varepsilon_{-1}^{(r+\ell)}, \varepsilon
_{-2}^{(r+\ell)}, \ldots.
\]
Consequently, the covariances
in (\ref{e:lr-cs}) vanish except when $r=k$ or $r+\ell=k$, so condition
(\ref{e:lr-cs}) always holds for linear processes.

For general nonlinear sequences, the difference
\[
X_k - X_{k}^{(k)} = f(\varepsilon_k, \ldots, \varepsilon_{1},
\varepsilon_{0},\varepsilon_{-1},
\ldots)
- f\bigl(\varepsilon_k, \ldots, \varepsilon_{1}, \varepsilon
_{0}^{(k)},\varepsilon_{-1}^{(k)}, \ldots\bigr),
\]
cannot be expressed only in terms of the errors (\ref{e:eg0}), but the
errors $\varepsilon_k, \ldots, \varepsilon_{1}$ should approximately
cancel, so that the
difference $X_k - X_{k}^{(k)}$ is small and very weakly correlated with
$X_r^{(r)} X_{r+\ell}^{(r+\ell)}$.

With this background, we now formulate the following result.
\begin{theorem} \label{t:lr-s}
Suppose $\{ X_n \}\in L^4$ is an $L^4$--$m$-approximable and assume condition
(\ref{e:lr-ct}) holds. If Assumption \ref{a:wq} holds,
then
$
\hat\sigma^2 \stackrel{P}{\rightarrow}\sum_{j=-\infty}^\infty
\gamma_j.
$
\end{theorem}

Theorem \ref{t:lr-s} is proven in Section \ref{s:p-si}. The general
plan of the proof is the same as that of the proof of Theorem 3.1 of
Giraitis et al. \cite{giraitiskokoszkaleipus2003RV}, but the
verification of
the crucial relation (\ref{e:lr3}) uses a new approach based on
$L^4$--$m$-approximability. The arguments preceding (\ref{e:lr3}) show
that replacing $\bar X_N$ by $\mu= EX_0$ does not change the limit.
We note that the condition $q^2/N \to0$ we assume is stronger than
the condition $q/N \to0$ assumed by
Giraitis et al. \cite{giraitiskokoszkaleipus2003RV}. This difference
is of
little practical consequence, as the optimal bandwidths for the
kernels used in practice are typically of the order $O(N^{1/5})$.
Finally, we notice that by further strengthening conditions on the
behavior of the bandwidth function $q=q(N)$, the convergence in
probability in Theorem \ref{t:lr-s} could be replaced by the almost
sure convergence, but we do not pursue this research here. The
corresponding result under condition (\ref{e:gklt-cond}) was established
by Berkes et al. \cite{berkeshorvathkokoszkashao2005}; it is also stated
without proof as part of Theorem A.1 of
Berkes et al. \cite{berkeshorvathkokoszkashao2006}.

We now turn to the vector case in which the data are of the
form
\[
\mathbf{X}_n = [X_{1n}, X_{2n}, \ldots, X_{dn}]^T,\qquad n=1,2,
\ldots, N.
\]
Just as in the scalar case, the estimation of the mean by the sample
mean does not affect the limit of the kernel long-run variance estimators,
so \textit{we assume that $E X_{in}=0$} and define the autocovariances
as
\[
\gamma_r(i,j) = E[X_{i0} X_{jr}],\qquad 1 \le i, j \le d.
\]
If $r\ge0$, $\gamma_r(i,j)$ is estimated by
$N^{-1} \sum_{n=1}^{N-r} X_{in}X_{j,n+r}$, but if $r<0$, it
is estimated by $N^{-1} \sum_{n=1}^{N-|r|} X_{i, n+|r|}X_{j,n}$.
We therefore define the autocovariance matrices
\[
\hat{\bolds\Gamma}_r=
\cases{
\displaystyle N^{-1} \sum_{n=1}^{N-r} \mathbf{X}_n \mathbf{X}_{n+r}^T, &\quad if
$r \ge0$,\cr
\displaystyle N^{-1} \sum_{n=1}^{N-|r|} \mathbf{X}_{n+ |r|} \mathbf{X}_{n}^T,
&\quad
if $r < 0$.}
\]
The variance
$\operatorname{Var}[N^{-1} \bar\mathbf{X}_n]$
has $(i,j)$-entry
\[
N^{-2} \sum_{m,n=1}^N E[X_{im} X_{jn}]
= N^{-1} \sum_{|r| <N} \biggl(1 - \frac{|r|}{N}\biggr)\gamma_r(i,j),
\]
so the long-run variance is
\[
{\bolds\Sigma}= \sum_{r=-\infty}^\infty{\bolds\Gamma}_r,\qquad
{\bolds\Gamma}_r := [\gamma_r(i,j), 1 \le i,j \le d],
\]
and its kernel estimator is
%
%e4.8 ###
%
\begin{equation}\label{e:kernelestimator}
\hat{\bolds\Sigma}= \sum_{|r| \le q} \omega_q(r) \hat{\bolds
\Gamma}_r.
\end{equation}

The consistency of $\hat{\bolds\Sigma}$ can be established by
following the lines
of the proof of Theorem \ref{t:lr-s} for every fixed entry of the
matrix $\hat{\bolds\Sigma}$. Conditition (\ref{e:lr-ct}) must be
replaced by
%
%e4.9 ###
%
\begin{equation}\label{e:lr-ctv}
N^{-1} \sum_{k,l=0}^{q(N)} \sum_{r=1}^{N-1}
\max_{1 \le i,j \le d} \bigl|\operatorname{Cov}\bigl( X_{i0} \bigl(X_{jk}- X_{jk}^{(k)}\bigr),
X_{ir}^{(r)} X_{j, r+\ell}^{(r+\ell)} \bigr)\bigr|\to0.
\end{equation}
Condition (\ref{e:lr-ctv}) is analogous to cumulant
conditions for vector processes which require summability
of fourth order cross-cumulants of all scalar components
(see, e.g., Andrews \cite{andrews1991}, Assumption A, page 823).

For ease of reference we state these results as a theorem.
\begin{theorem} \label{t:lr-v}
\textup{(a)} If $\{ \mathbf{X}_n \}\in L_{\mathbb{R}^d}^2$ is an
$L^2$--$m$-approximable
sequence, then the series $\sum_{r=-\infty}^\infty{\bolds\Gamma}_r$
converges absolutely.
\textup{(b)} Suppose $\{ \mathbf{X}_n \}\in L_{\mathbb{R}^d}^4$ an
$L^4$--$m$-approximable
sequence such that condition
(\ref{e:lr-ctv}) holds. If Assumption \ref{a:wq} holds,
then $\hat{\bolds\Sigma}\stackrel{P}{\rightarrow}{\bolds\Sigma}$.
\end{theorem}

We are now able to turn to
functional data. Suppose $\{X_n\}\in L_H^2$ is a zero mean sequence,
and $v_1, v_2, \ldots, v_d$ is any set of orthonormal functions in $H$.
Define $X_{in} = \int X_n(t) v_i(t)\,dt$,
$\mathbf{X}_n =[X_{1n}, X_{2n}, \ldots, X_{dn}]^T$ and
$\mathbf{\Gamma}_r=\operatorname{Cov}(\mathbf{X}_0,\mathbf{X}_r)$.
A direct verification shows that if $\{X_n\}$ is $L^p$--$m$-approximable,
then so is the vector sequence $\{\mathbf{X}_n \}$. We thus obtain the
following
corollary.
\begin{corollary} \label{c:lr-f}
\textup{(a)} If $\{ X_n \}\in L_H^2$ is an $L^2$--$m$-approximable sequence,
then the series $\sum_{r=-\infty}^\infty{\bolds\Gamma}_r$
converges absolutely.
\textup{(b)} If, in addition,
$\{ X_n \}$ is $L^4$--$m$-approximable and Assumption \ref{a:wq}
and condition (\ref{e:lr-ctv}) hold, then
$\hat{\bolds\Sigma}\stackrel{P}{\rightarrow}{\bolds\Sigma}$.
\end{corollary}

In Corollary \ref{c:lr-f},
the functions $v_1, v_2, \ldots, v_d$
form an arbitrary orthonormal deterministic basis. In many
applications, a random basis consisting of the
estimated principal components ${\hat v}_1, {\hat v}_2, \ldots, {\hat
v}_d$ is used.
The scores with respect to this basis are defined by
\[
\hat{\eta}_{\ell i}=\int\bigl(X_i(t)-\bar{X}_N(t)\bigr){\hat v}_\ell(t)
\,dt,\qquad 1\leq \ell\leq d.
\]
To use the results established so far,
it is convenient to decompose the stationary sequence $\{ X_n \}$ into
its mean and a zero mean process; that is, we set $X_n(t)= \mu(t) + Y_n(t)$,
where $EY_n(t)=0$. We introduce
the unobservable quantities
%
%e4.10 ###
%
\begin{equation}\label{e:beta}\qquad
\beta_{\ell n }=\int Y_n(t) v_\ell(t) \,dt,\qquad
\hat\beta_{\ell n}=\int Y_n(t) {\hat v}_\ell(t) \,dt,\qquad
1\leq\ell\leq d.
\end{equation}
We then have the following proposition which will be useful in
most statistical procedures for functional time series. An application
to change point detection is developed in Section \ref{s:change}.
\begin{proposition}\label{p:diffbartlett}
Let $\hat{\mathbf{C}}=\operatorname{diag}(\hat{c}_1,\ldots,\hat{c}_d)$,
with $\hat{c}_i= \operatorname{sign}(\langle v_i, {\hat v}_i \rangle)$.
Suppose $\{ X_n \}\in L_H^4$ is $L^4$--$m$-approximable and that
(\ref{e:la>d}) holds. Assume further that Assumption \ref{a:wq}
holds with a stronger condition $q^4/N\to0$. Then
\[
|\hat{{\bolds\Sigma}}(\bolds\beta)-\hat{{\bolds\Sigma}}(\hat
{\mathbf{C}}\hat{\bolds\beta}
)|=o_P(1)
\quad\mbox{and}\quad
|\hat{{\bolds\Sigma}}(\hat{{\bolds\eta}})-\hat{{\bolds\Sigma
}}(\hat{\bolds\beta})|=o_P(1).
\]
\end{proposition}

The proof of Proposition \ref{p:diffbartlett} is delicate and is
presented in Section \ref{s:p-si}. We note that condition
(\ref{e:lr-ctv}) does not appear in the statement of Proposition
\ref{p:diffbartlett}. Its point is that if $\hat{{\bolds\Sigma
}}(\bolds\beta)$ is
consistent under some conditions, then so is $\hat{{\bolds\Sigma
}}(\hat{{\bolds\eta}})$.

%Let us finally note that
%working out a theory of long-run variance estimation for functional
%data, without any projections, would not be out of reach, the
%derivations
%would be similar to the scalar case, with functions and operators in
%place
%of the scalars $X_n$ and $\sigma^2$.
%Such results are, however, not immediately
%useful for statistical applications. In most statistical procedures in
%current use, the dimension of the functional data is initially
%reduced in one way
%or another by projecting on fixed or random basis functions,
%which optimally describe the data in a sense which depends on an
%application.
%}

%s5 ###
\section{Change point detection}\label{s:change}

Functional time series are obtained from data collected
sequentially over time, and it is natural to expect that conditions
under which observations are made may change. If this is the case,
procedures developed for stationary series will produce spurious
results. In this section, we develop a procedure for the detection
of a change in the mean function of a functional time series, the
most important possible change.
In addition to its practical relevance, the requisite theory illustrates
the application of the results developed in Sections~\ref{ss:conv} and
\ref{s:lr}. The main results of this Section, Theorems \ref{t:QFFCLT}
and \ref{th:alternative}, are proven in Section \ref{s:pQ}.

We thus consider testing the null hypothesis,
\[
H_0\dvtx E X_1(t) = E X_2(t) = \cdots= E X_N(t),\qquad t \in[0,1].
\]
Note that under $H_0$, we do not specify the value of the common mean.

Under the alternative, $H_0$ does not hold. The test we construct
has a particularly good power against the alternative in which the
data can be divided into several consecutive segments, and the
mean is constant within each segment but changes from segment to segment.
The simplest case of only two segments (one change point) is specified
in Assumption \ref{a:con-si}.
First we note that
under the null hypothesis, we can represent each functional observation
as
%
%e5.1 ###
%
\begin{equation}\label{e:Xdec}
X_i(t)= \mu(t) + Y_i(t),\qquad E Y_i(t) =0.
\end{equation}
The following assumption specifies conditions on $\mu(\cdot)$ and the
errors $Y_i(\cdot)$ needed to establish the convergence of the test
statistic under $H_0$.
\begin{assumption}\label{a:null-si}
The mean $\mu$ in (\ref{e:Xdec}) is in $H$.
The error functions $Y_i\in L_H^4$ are
$L^4$--$m$-approximable mean zero random elements
such that the eigenvalues of their covariance operator
satisfy (\ref{e:la>d}).
\end{assumption}

Recall that the $L^4$--$m$-approximability implies that the
$Y_i$ are identically distributed with
$\nu_4(Y_i)<\infty$. In particular, their covariance function,
\[
c(t,s) =E [Y_i(t) Y_i(s)],\qquad 0 \le t,s\le1,
\]
is square integrable, that is, is in $L^2([0,1]\times[0,1])$.

We develop the theory under the alternative of exactly one change point,
but the procedure is applicable to multiple change points by using
a segmentation algorithm described in
Berkes et al. \cite{berkesgabryshorvathkokoszka2009} and dating back
at least to
Vostrikova \cite{vostr1981}.
\begin{assumption} \label{a:con-si}
The observations follow the model
\[
X_i(t) =
\cases{\mu_1(t) + Y_i(t), &\quad $1 \le i\le k^*$,\cr
\mu_2(t) + Y_i(t), &\quad $k^* < i\le N$,}
\]
in which the $Y_i$ satisfy Assumption \ref{a:null-si}, the mean
functions $\mu_1$ and $\mu_2$ are in $L^2$ and
\[
k^*=[n\theta] \qquad\mbox{for some } 0<\theta<1.
\]
\end{assumption}

The general idea of testing is similar to that
developed in Berkes et al. \cite{berkesgabryshorvathkokoszka2009}
for independent observations, the central difficulty is in accommodating
the dependence. To define the
test statistic, recall that bold symbols denote $d$-dimensional
vectors, for example,
$\hat{{\bolds\eta}}_i=[\hat{\eta}_{1i},\hat{\eta}_{2i},\ldots
,\hat{\eta}_{di}]^T$.
To lighten the notation, define the partial sums process,
$\mathbf{S}_N(x,\bolds\xi)=\sum_{n=1}^{\lfloor Nx\rfloor}\bolds
\xi_n$, $x\in[0,1]$,
and the process, $\mathbf{L}_N(x,\bolds\xi)=\mathbf{S}_N(x,\bolds
\xi)-x\mathbf{S}_N(1,\bolds\xi)$,
where $\{\bolds\xi_n \}$ is a generic $R^d$-valued sequence.
Denote by ${\bolds\Sigma}(\bolds\xi)$ the long-run variance of the
sequence $\{\bolds\xi_n \}$, and by $\hat{{\bolds\Sigma}}(\bolds
\xi)$
its kernel estimator (see Section \ref{s:lr}).
The proposed test statistic is then
%
%e5.2 ###
%
\begin{equation}\label{e:test}
T_N(d)=\frac{1}{N}\int_0^1\mathbf{L}_N(x,\hat{{\bolds\eta}})^T
\hat{{\bolds\Sigma}}(\hat{{\bolds\eta}})^{-1}
\mathbf{L}_N(x,\hat{{\bolds\eta}}) \,dx.
\end{equation}
Our first theorem establishes its asymptotic null distribution.
\begin{theorem}\label{t:QFFCLT} Suppose $H_0$ and
Assumption \ref{a:null-si} hold.
If the estimator $\hat{{\bolds\Sigma}}(\hat{{\bolds\eta}})$ is
consistent, then
%
%e5.3 ###
%
\begin{equation}\label{teststat}
T_N(d)\stackrel{d}{\rightarrow}T(d):=
\sum_{\ell=1}^d\int_0^1B_\ell^2(x) \,dx,
\end{equation}
where
$\{B_\ell(x), x\in[0,1]\}$, $1\leq\ell\leq d$ are independent
Brownian bridges.
\end{theorem}

The distribution of the random variable $T(d)$ was derived
by Kiefer \cite{kiefer1959}.
%The distribution function is complicated,
%and therefore we give for the convenience of the reader
%the critical values when $d=1,\ldots,5$ in Table \ref{tab:critvalues}.
%It is known that $ET(d)=d/6$ and $\mathrm{Var}(T(d))=d/45$, so a normal
%approximation can be used, Aue et al.
%show that it is applicable for $d\ge10$.
%Simulated finite sample
%critical values for small $d$'s are tabulated in
%Berkes et al. \cite{berkesgabryshorvathkokoszka2009}.
%
%levels $1-\alpha$.}
% \hline
% % after \ \hline or \cline{col1-col2} \cline{col3-col4}...
% $\alpha\setminus d$ & 1 & 2 & 3 & 4 & 5 \\
% \hline
% 0.1 &.35 & .61 & .84 & 1.06 & 1.28 \\
% 0.05 & .46 & .75 & 1.00 & 1.24 & 1.47 \\
% 0.01 & .74 & 1.07 & 1.36 & 1.62 & 1.88 \\
% \hline
%}%comment
The limit distribution is the same as in the case of independent observations;
this is possible because the long-run variance estimator
$\hat{{\bolds\Sigma}}(\hat{{\bolds\eta}})$ soaks up the dependence.
Sufficient conditions for its consistency are stated
in Section \ref{s:lr}, and, in addition to the assumptions
of Theorem \ref{t:QFFCLT}, they are: Assumption \ref{a:wq} with
$q^4/N\to0$,
and condition (\ref{e:lr-ctv}).

The next result shows that our test has asymptotic power 1.
Our proof requires the following condition:
%
%e5.4 ###
%
\begin{equation}\label{e:conalt}
\hat{\bolds\Sigma}(\hat{\bolds\eta})
\stackrel{\mathrm{a.s.}}{\rightarrow}\bolds\Omega\qquad
\mbox{where }
\bolds\Omega\mbox{ is some
positive definite matrix.}
\end{equation}

Condition (\ref{e:conalt}) could be replaced by weaker technical conditions,
but we prefer it, as it
leads to a transparent, short proof. Essentially, it states that
the matrix $\hat{\bolds\Sigma}(\hat{\bolds\eta})$ does
not become
degenerate in the limit, and the matrix $\bolds\Omega$ has only
positive eigenvalues. A condition like
(\ref{e:conalt}) is not needed for independent $Y_i$
because that case does not require normalization with
the long-run covariance matrix.
To formulate our result, introduce
vectors $\bolds\mu_1, \bolds\mu_2 \in\mathbb{R}^d$
with coordinates
\[
\int\mu_1(t){v}_\ell(t)\,dt \quad\mbox{and}\quad
\int\mu_2(t){v}_\ell(t)\,dt,\qquad 1\leq\ell\leq d.
\]
\begin{theorem}\label{th:alternative} Suppose Assumption \ref{a:con-si}
and condition (\ref{e:conalt}) hold. If the vectors $\bolds\mu_1$
and $\bolds\mu_2$
are not equal, then
$
T_N(d)\stackrel{P}{\rightarrow}\infty.
$
\end{theorem}

We conclude this section with two numerical examples which illustrate
the effect of dependence on our change point detection procedure.
Example \ref{ex:simlongrun} uses synthetic data while
Example \ref{ex:pm10} focuses on particulate pollution data.
Both show that using statistic (\ref{e:test})
with $\hat{\bolds{\Sigma}}(\hat{\bolds{\eta}})$ being
the estimate for just
the covariance, not the long-run
covariance matrix, leads to spurious rejections
of $H_0$, a nonexistent change point can be detected with
a large probability.
\begin{example}\label{ex:simlongrun}
We simulate 200 observations of the functional AR(1) process
of Example \ref{ex:m-far}, when $\Psi$ has the parabolic integral kernel
$\psi(t,s)=\gamma\cdot(2-(2x-1)^2-(2y-1)^2 )$.
We chose the constant $\gamma$ such that $\|\Psi\|_{\mathcal{S}}=0.6$
(the Hilbert--Schmidt norm).
The innovations $\{\varepsilon_n\}$ are standard Brownian
bridges. The first 3 principal components
explain approximately 85\% of the total variance,
so we compute the test statistic $T_{200}(3)$ given in (\ref{e:test}).
For the estimation of the long-run covariance matrix
$\bolds{\Sigma}$ we use the {Bartlett kernel}
\[
\omega_q^{(1)}(j)=\cases{
1-|j|/(1+q), &\quad if $|j|\leq q$; \cr
0, &\quad otherwise.}
\]
We first let $q=0$ which corresponds to using just
the sample covariance of $\{\bolds{\hat{\eta}}_n\}$
in the normalization for the test statistic (\ref{e:test})
(dependence is ignored).
We use 1000 replications and the 5\% confidence level.
The rejection rate is $23.9\%$,
much higher than the nominal level of $5\%$.
In contrast, using
an appropriate estimate for the long-run variance,
the reliability of the test improves dramatically.
Choosing an optimal bandwidth $q$ is a separate problem
which we do not pursue here. Here we
adapt the formula
$q\approx1.1447 (a N)^{1/3}$, $a = \frac{4\psi^2}{(1+\psi)^4}$
valid for a a scalar AR(1) process with the autoregressive
%
%f2 ###
%
\begin{figure}[b]\vspace*{-4pt}

\includegraphics{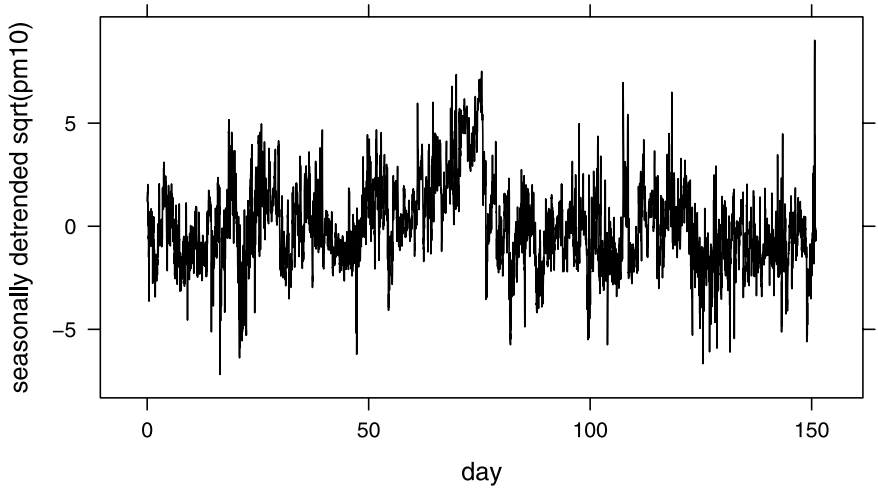}

\caption{Seasonally detrended $\sqrt{\mathtt{pm10}}$,
Nov 1, 2008--Mar 31, 2009.}
% Lower panels: Original intraday curve (left) and
%smoothed intraday curve (rigth).}
\label{fig:pm10}
\end{figure}
coefficient $\psi$ (Andrews \cite{andrews1991}). Using this formula with
$\psi= \|\Psi\|_{\mathcal{S}}=0.6$ results in $q=4$. This choice
gives the empirical rejection rate of $3.7\%$, much closer to the
nominal rate of $5\%$.
\end{example}
\begin{example}\label{ex:pm10}
This example, which uses \texttt{pm10} (particulate matter with
diameter $<10$ $\mu\mathrm{m}$, measured in $\mu\mathrm{g}/\mathrm{m}^3$)
data, illustrates a similar phenomenon as
Example \ref{ex:simlongrun}. For the analysis we use \texttt{pm10}
concentration data measured in the Austrian city of Graz during the
winter of 2008/2009 ($N$=151). The data are given in 30 minutes
resolution, yielding an intraday frequency of 48 observations. As
in Stadtlober, H{\"o}rmann and Pfeiler \cite
{stadloberhoermannpfeiler2008} we use a square root
transformation to reduce heavy tails. Next we remove possible
weekly periodicity by subtracting the corresponding mean vectors
obtained from the different weekdays. A time series plot of this
new sequence is given in Figure \ref{fig:pm10}. The data look
relatively stable, although a shift appears to be possible in the
center of the time series. It should be emphasized, however, that
\texttt{pm10} data, like many geophysical time series, exhibit a
strong, persistent, positive autocorrelation structure. These series
are stationary over long periods of time with an appearance of
local trends or shifts at various time scales (random self-similar
or fractal structure).

The daily measurement vectors are transformed into smooth functional
data using 15 B-splines functions of order 4. The functional
principal component analysis yields that the first three principal
components explain $\approx84\%$ of the total variability, so we
use statistic (\ref{e:test}) with $d=3$. A look at the \texttt{acf} and
\texttt{pacf} of the first empirical PC scores
(Figure \ref{fig:pm10acf}) suggests an AR(1), maybe AR(3) behavior.
The second and third empirical PC scores show no significant
autocorrelation structure. We use the formula given in
Example \ref{ex:simlongrun} with $\psi=0.70$ (\texttt{acf} at lag 1)
and $N=151$ and obtain $q\approx4$. This gives $T_{151}(3)=0.94$ which
%
%f3 ###
%
\begin{figure}

\includegraphics{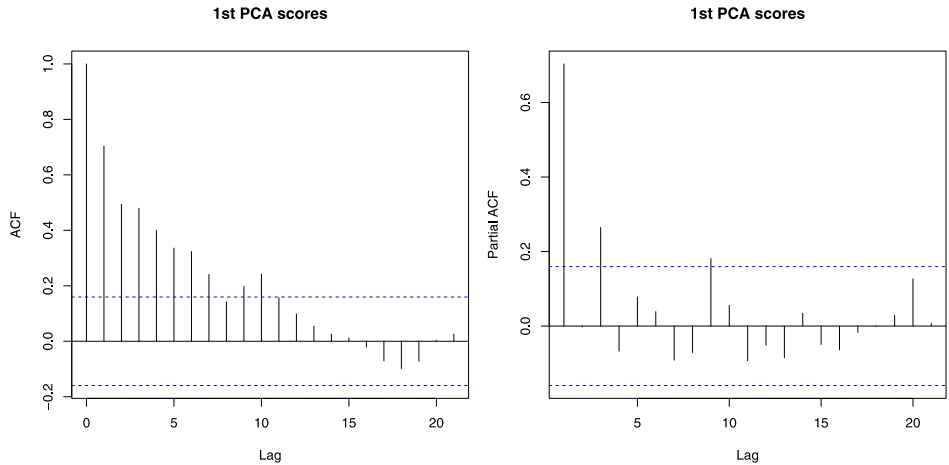}

\caption{Left panel: sample autocorrelation function
of the first empirical PC scores. Right panel: sample partial autocorrelation
function of the first empirical PC scores.}\label{fig:pm10acf}
\end{figure}
is close to the critical value $1.00$ when testing at a $95\%$
confidence level but does not support rejection of the no-change
hypothesis. In contrast, using only the sample covariance matrix in
(\ref{teststat}) gives $T_{151}(3)=1.89$ and thus a clear and
possibly wrongful rejection of the null hypothesis.
\end{example}

%s6 ###
\section{Functional linear model with
dependent regressors} \label{s:flm-si}

The functional linear model is one of the most widely used tools
of FDA. Its various
forms are introduced in Chapters 12--17
of Ramsay and Silverman \cite{ramsaysilverman2005}. To name a few
recent references
we mention
Cuevas, Febrero and Fraiman \cite{cuevasfebrerofreiman2002},
Malfait and Ramsay~\cite{malfaitramsay2003},
Cardot et al. \cite{cardotfms2003},
Cardot, Ferraty and Sarda \cite{cardotferratysarda2003},
Chiou, M{\"u}ller and Wang \cite{chioumullerwang2004},
M{\"u}ller and Stadtm{\"u}ller \cite{mullerstadtmuller2005},
Yao, M{\" u}ller and Wang \cite{yaomullerwang2005},
Cai and Hall~\cite{caihall2006},
Chiou and M{\"u}ller \cite{chioumuller2007}, Li and Hsing \cite{lihsing2007},
Reiss and Ogden \cite{reissogden2007}, Reiss and Ogden
\mbox{\cite{reissogden2009a,reissogden2009b}}.

We focus on the fully functional model of the form
%
%e6.1 ###
%
\begin{equation}\label{e:fm-i}
Y_n(t) = \int\psi(t,s) X_n(s) + \varepsilon_n(t),\qquad n=1,2
\ldots, N,
\end{equation}
in which both the regressors and the responses are functions.
The results of this section can be easily specialized to the
case of scalar responses.

In (\ref{e:fm-i}), the regressors are random functions,
assumed to be independent and identically distributed. As
explained in Section \ref{s:sie}, for functional time series
the assumption of the independence of the $X_n$ is often
questionable, so it is important to investigate if procedures
developed and theoretically justified for independent regressors
can still be used if the regressors are dependent.

We focus here on the estimation of the kernel $\psi(t,s)$. Our result
is motivated by the work of Yao, M{\" u}ller and Wang \cite
{yaomullerwang2005} who
considered functional regressors and responses obtained from sparce
\textit{independent} data measured with error. The data that motivates
our work are measurements of physical quantities obtained with
negligible errors or financial transaction data obtained without
error. In both cases the data are available at fine time grids, and
the main concern is the presence of temporal dependence between the
curves $X_n$. We therefore merely assume that the sequence $\{X_n\}\in L_H^4$
is $L^4$--$m$-approximable, which, as can be easily seen, implies the
$L^4$--$m$-approximability of $\{Y_n\}$.
To formulate additional
technical assumptions, we need to introduce some notation.

We assume that the errors $\varepsilon_n$ are i.i.d. and independent
of the
$X_n$, and denote by $X$ and $Y$ random functions with the same
distribution as $X_n$ and $Y_n$, respectively.
We work with their expansions
\[
X(s) = \sum_{i=1}^\infty\xi_i v_i(s),\qquad
Y(t) = \sum_{j=1}^\infty\zeta_j u_j(t),
\]
where the $v_j$ are the FPCs of $X$ and the $u_j$ the FPCs of $Y$,
and
$
\xi_i = \langle X, v_i \rangle,
\zeta_j = \langle Y, u_j \rangle.
$
Indicating with the ``hat'' the corresponding empirical quantities,
an estimator of $\psi(t,s)$ proposed by Yao, M{\" u}ller and Wang
\cite{yaomullerwang2005} is
\[
\hat\psi_{K L}(t,s) = \sum_{k=1}^K \sum_{\ell= 1}^L
{\hat\lambda}_\ell^{-1} \hat\sigma_{\ell k} \hat u_k(t) {\hat
v}_\ell(s),
\]
where $\hat\sigma_{\ell k}$ is an estimator of $E[\xi_\ell\zeta_k]$.
We will work with the simplest estimator,
%
%e6.2 ###
%
\begin{equation}\label{e:hat-sigma-lk}
\hat\sigma_{\ell k} = \frac{1}{N} \sum_{i=1}^N \langle X_i,
{\hat v}_\ell
\rangle
\langle Y_i, \hat u_k \rangle,
\end{equation}
but any estimator for which Lemma \ref{l:consistentsigma}
holds can be used without affecting the rates.

Let $\lambda_j$ and $\gamma_j$ be the eigenvalues corresponding to
$v_j$ and $u_j$.
Define $\alpha_j$ as in Lemma \ref{l:B4.3}, and define $\alpha_j'$
accordingly with
$\gamma_j$ instead of $\lambda_j$.
Set
\[
h_L=\min\{\alpha_j, 1\leq j\leq L\},\qquad
h_L'=\min\{\alpha_j', 1\leq j\leq L\}.
\]
To establish the consistency of the estimator $\hat\psi_{K L}(t,s)$
we assume that
%
%e6.3 ###
%
\begin{equation}\label{e:ymw}
\Psi:=\sum_{k=1}^\infty\sum_{\ell=1}^\infty
\frac{(E[\xi_{\ell} \zeta_{k}])^2}{\lambda_\ell^2} < \infty
\end{equation}
and that the following assumption holds:
\begin{assumption}\label{a:linreg}
(i) We have $\lambda_1>\lambda_2>\cdots$ and $\gamma_1>\gamma
_2>\cdots.$\vspace*{1pt}

(ii) We have $K=K(N)$, $L=L(N)\to\infty$ and
$
\frac{KL}{\lambda_L\min\{h_K,h_L'\}} =o (N^{1/2} ).
$
\end{assumption}

For model (\ref{e:fm-i}), condition (\ref{e:ymw}) is equivalent
to the assumption that $\psi(t,s)$ is a Hilbert--Schmidt kernel,
that is, $\dint\psi^2(t,s) \,dt\,ds< \infty$. It is formulated in
the same way
as in Yao, M{\" u}ller and Wang \cite{yaomullerwang2005} because this
form is convenient in the
theoretical arguments. Assumption \ref{a:linreg} is much shorter than
the corresponding assumptions of Yao, M{\" u}ller and Wang \cite
{yaomullerwang2005} which take
up over two pages. This is because we do not deal with smoothing and so can
isolate the impact of the magnitude of the eigenvalues on the bandwidths
$K$ and $L$.
\begin{theorem}\label{t:linregdep}
Suppose $\{ X_n \}\in L_H^4$ is a zero mean $L^4$--$m$-approximable sequence
independent of the sequence of i.i.d. errors $\{\varepsilon_n\}$.
If (\ref{e:ymw}) and Assumption~\ref{a:linreg} hold, then
%
%e6.4 ###
%
\begin{equation}\label{e:ymw-con}
\dint[\hat\psi_{K L}(t,s) - \psi(t,s) ]^2 \,dt\,ds \stackrel
{P}{\rightarrow}0,\qquad
(N \to\infty).
\end{equation}
\end{theorem}

The proposition of Theorem \ref{t:linregdep} is comparable to the first
part of Theorem 1 in Yao, M{\" u}ller and Wang \cite{yaomullerwang2005}.
Both theorems are established under (\ref{e:ymw}) and
finite fourth moment conditions. Otherwise the settings
are quite different. Yao, M{\" u}ller and Wang \cite{yaomullerwang2005}
work under the assumption
that the subject $(Y_i,X_i)$, $i=1,2,\ldots$ are independent and sparsely
observed
whereas the crucial point of our approach is that we
allow dependence. Thus Theorems 1 and 2 in the related
paper
Yao, M{\" u}ller and Wang \cite{yaomullerwang2005a}, which serve as the basic
ingredients for their results, cannot be used here
and have to be replaced directly
with the theory developed in Section \ref{ss:conv}
of this paper. Furthermore, our proof goes without
complicated assumptions on the resolvents of
the covariance operator, in particular without
the very technical assumptions (B.5) of Yao, M{\" u}ller and Wang \cite
{yaomullerwang2005}.
In this sense, our short alternative proof might be of value even in the
case of independent observations.

\begin{appendix}\label{app}
\section*{Appendix}

We present the proofs of results stated in
Sections \ref{ss:conv}--\ref{s:flm-si}. Throughout we
will agree on the following conventions. All $X_n\in L_H^2$
satisfy Assumption \ref{a:reg}. A~generic $X$, which
is assumed to be equal in distribution
to $X_1$, will be used
at some places. Any constants occurring will be denoted by
$\kappa_1, \kappa_2,\ldots.$ The $\kappa_i$ may
change their values from proof to proof.

%s6.1 ###
\subsection{\texorpdfstring{Proofs of the results of Sections \protect\ref{ss:conv} and
\protect\ref{s:lr}}{Proofs of the results of Sections 3 and 4}}\label{s:p-si}

\mbox{}

\begin{pf*}{Proof of Theorem \protect\ref{t:ChatW}}
We assume for simplicity that $EX=0$ and set
\[
{\hat c}(t,s) = N^{-1} \sum_{n=1}^N X_n(t) X_n(s),\qquad
c(t,s) = E[X(t) X(s)].
\]
The proof with a general mean function $\mu(t)$ requires some additional
but similar arguments.
The Cauchy--Schwarz inequality shows that ${\hat c}(\cdot,\cdot)$ and
$c(\cdot,\cdot)$ are Hilbert--Schmidt
kernels, so $\hat C - C$ is a Hilbert--Schmidt operator with
the kernel ${\hat c}(t,s) - c(t,s)$. Consequently,
\[
NE\|\hat C - C\|_{{\mathcal S}}^2
= N \dint\operatorname{Var}
\Biggl[ N^{-1} \sum_{n=1}^N \bigl(X_n(t) X_n(s) - E[X_n(t) X_n(s)]\bigr) \Biggr] \,dt\,ds.
\]
For fixed $s$ and $t$, set
$
Y_n = X_n(t) X_n(s) - E[X_n(t) X_n(s)].
$
Due the stationarity of the sequence $\{ Y_n \}$ we have
\[
\operatorname{Var}\Biggl( N^{-1} \sum_{n=1}^N Y_n\Biggr)= N^{-1} \sum_{|r|<N}
\biggl(1 - \frac{|r|}{N}\biggr)\operatorname{Cov}(Y_1, Y_{1+r})
\]
and so
\[
N\operatorname{Var}\Biggl( N^{-1} \sum_{n=1}^N Y_n\Biggr)
\le\operatorname{Var}(Y_1) +{ 2 \sum_{r=1}^\infty}|{\operatorname
{Cov}}(Y_1, Y_{1+r})|.
\]
Setting
$
Y_n^{(m)} = X_n^{(m)}(t) X_n^{(m)}(s) - E[X_n(t) X_n(s)],
$
we obtain
\[
|{\operatorname{Cov}(Y_1, Y_{1+r})}| = \bigl|{\operatorname{Cov}}\bigl(Y_1,
Y_{1+r}-Y_{1+r}^{(r)}\bigr)\bigr|
\le[\operatorname{Var}(Y_1)]^{1/2} \bigl[\operatorname
{Var}\bigl(Y_{1+r}-Y_{1+r}^{(r)}\bigr)\bigr]^{1/2}.
\]
Consequently,
$
NE\|\hat C - C\|_{{\mathcal S}}^2
$
is bounded from above by
\begin{eqnarray*}
&&\dint\operatorname{Var}[X(t)X(s)]\,dt\,ds\\
&&\qquad{}+ 2 \sum_{r=1}^\infty
\dint[\operatorname{Var}
(X(t)X(s))]^{1/2}
\\
&&\qquad\hspace*{51.5pt}{} \times\bigl[\operatorname{Var}\bigl(X_{1+r}(t)X_{1+r}(s)-X_{1+r}^{(r)}(t)
X_{1+r}^{(r)}(s) \bigr)\bigr]^{1/2}\,dt\,ds.
\end{eqnarray*}
For the first summand we have the upper bound $\nu_4^4(X)$ because
%
%e6.5 ###
%
\setcounter{equation}{0}
\begin{equation}\label{e:si1}
\dint E[X^2(t) X^2(s)]\,dt\,ds =
E \int X^2(t)\,dt \int X^2(s)\,ds = \nu_4^4(X).
\end{equation}
To find upper bounds for the summands
in the infinite sum, we use the inequality
%
%e6.6 ###
%
\begin{equation}\label{e:abcd}
|ab-cd|^2\leq2a^2(b-d)^2+2d^2(a-c)^2,
\end{equation}
which yields
\begin{eqnarray*}
&& \dint[\operatorname{Var}(X(t)X(s))]^{1/2}
\bigl[\operatorname{Var}\bigl(X_{1+r}(t)X_{1+r}(s)-X_{1+r}^{(r)}(t) X_{1+r}^{(r)}(s)
\bigr)\bigr]^{1/2}\,dt\,ds\\
&&\qquad\le\dint[E(X^2(t)X^2(s))]^{1/2}
\bigl[E\bigl( X_{1+r}(t)X_{1+r}(s)\\
&&\hspace*{153.73pt}{} -X_{1+r}^{(r)}(t) X_{1+r}^{(r)}(s)
\bigr)^2\bigr]^{1/2}\,dt\,ds\\
&&\qquad\le\sqrt{2}\dint[E(X^2(t)X^2(s))]^{1/2}
\bigl[EX_{1+r}^2(t)\bigl(X_{1+r}(s)- X_{1+r}^{(r)}(s)\bigr)^2\bigr]^{1/2}\,dt\,ds\\
&&\qquad\quad{} + \sqrt{2}\dint[E(X^2(t)X^2(s))]^{1/2}\\
&&\hspace*{77.3pt}{}\times
\bigl[EX_{1+r}^{(r)2}(s)\bigl(X_{1+r}(t)
-X_{1+r}^{(r)}(t)\bigr)^2\bigr]^{1/2}\,dt\,ds.
\end{eqnarray*}

For the first term, using the Cauchy--Schwarz inequality
and (\ref{e:si1}), we obtain
\begin{eqnarray*}
&&\dint[E(X^2(t)X^2(s))]^{1/2}
\bigl[EX_{1+r}^2(t)\bigl(X_{1+r}(s)- X_{1+r}^{(r)}(s)\bigr)^2\bigr]^{1/2}\,dt\,ds\\
&&\qquad \le\nu_4^2(X)
\biggl\{ E \biggl[\int X_{1+r}^2(t) \,dt \int\bigl(X_{1+r}(s)- X_{1+r}^{(r)}(s)\bigr)^2
\,ds \biggr]
\biggr\}^{1/2}\\
&&\qquad \le\nu_4^2(X) \{ E\|X_{1+r}\|^4\}^{1/4}
\bigl\{ E \bigl\|X_{1+r}(s)- X_{1+r}^{(r)}(s)\bigr\| \bigr\}^{1/4}\\
&&\qquad = \nu_4^3(X) \nu_4 \bigl(X_{1}- X_{1}^{(r)}\bigr).
\end{eqnarray*}
The exact same argument applies for the second term.
The above bounds imply (\ref{e:s2}).
\end{pf*}
\begin{pf*}{Proof of Theorem \protect\ref{t:lr-s}}
As in Giraitis et al. \cite{giraitiskokoszkaleipus2003RV}, set
$\mu= E X_0$ and
\begin{eqnarray*}
\tilde\gamma_j
&=& \frac{1}{N}\sum_{i=1}^{N-|j|} (X_i - \mu) \bigl(X_{i+|j|} - \mu\bigr),
\\
S_{k, \ell} &=& \sum_{i=k}^\ell(X_i - \mu).
\end{eqnarray*}
Observe that
\[
\hat\gamma_j - \tilde\gamma_j
= \biggl(1 - \frac{|j|}{N} \biggr) (\bar X_N - \mu)^2
+ \frac{1}{N} (\bar X_N - \mu) \bigl( S_{1, N-|j|} + S_{|j|+1, N} \bigr) =:
\delta_j.
\]
We therefore have the decomposition
\[
\hat\sigma^2 = \sum_{|j| \le q} \omega_q(j) \tilde\gamma_j
+ \sum_{|j| \le q} \omega_q(j) \delta_j =: \hat\sigma_1^2 + \hat
\sigma_2^2.
\]
The proof will be complete once we have shown that
%
%e6.7 ###
%
\begin{equation}\label{e:lr1}
\hat\sigma_1^2 \stackrel{P}{\rightarrow}\sum_{j=-\infty}^\infty
\gamma_j
\end{equation}
and
%
%e6.8 ###
%
\begin{equation}\label{e:lr2}
\hat\sigma_2^2 \stackrel{P}{\rightarrow}0.
\end{equation}

We begin with the verification of the easier relation
(\ref{e:lr2}). By (\ref{e:og-b}),
\begin{eqnarray*}
E|\hat\sigma_2^2 | &\leq& b \sum_{|j| \le q} E |\delta_j|\\
&\le& b \sum_{|j| \le q} E(\bar X_N - \mu)^2\\
&&{} + \frac{b}{N} [ E (\bar X_N - \mu)^2]^{1/2}
\sum_{|j| \le q} \bigl[ E\bigl( S_{1, N-|j|} + S_{|j|+1, N} \bigr)^2\bigr]^{1/2}.
\end{eqnarray*}
By Lemma \ref{l:cov-sum},
\[
E(\bar X_N - \mu)^2 = \frac{1}{N} \sum_{|j| \le N}
\biggl(1 - \frac{|j|}{N} \biggr)\gamma_j = O(N^{-1}).
\]
Similarly
$
E( S_{1, N-|j|} + S_{|j|+1, N} )^2 = O(N).
$
Therefore,
\[
E|\hat\sigma_2^2 | = O(qN^{-1}+ N^{-1}N^{-1/2} q N^{1/2}) = O(q/N).
\]

We now turn to the verification of (\ref{e:lr1}).
We will show that $E\hat\sigma_1^2 \to\sum_{j} \gamma_j$
and $\operatorname{Var}[\hat\sigma_1^2]\to0$.

By (\ref{e:og-c}),
\[
E \hat\sigma_1^2 = \sum_{|j| \le q} \omega_q(j) \frac{N-|j|}{N}
\gamma_j
\to\sum_{j= -\infty}^\infty\gamma_j.
\]

By (\ref{e:og-b}), it remains to show that
%
%e6.9 ###
%
\begin{equation}\label{e:lr3}
{\sum_{|k|, |\ell| \le q}} |{ \operatorname{Cov}}(\tilde\gamma_k,
\tilde\gamma_\ell)|
\to0.
\end{equation}
To lighten the notation, without any loss of generality,
\textit{we assume from now on that $\mu=0$}, so that
\[
\operatorname{Cov}(\tilde\gamma_k, \tilde\gamma_\ell) = \frac{1}{N^2}
\operatorname{Cov}\Biggl(\sum_{i=1}^{N-|k|}X_i X_{i+|k|},
\sum_{j=1}^{N-|\ell|}X_j X_{j+|\ell|} \Biggr).
\]
Therefore, by stationarity,
\begin{eqnarray*}
|{\operatorname{Cov}}(\tilde\gamma_k, \tilde\gamma_\ell)|
&\le&\frac{1}{N^2} \sum_{i, j=1}^N
\bigl|\operatorname{Cov}\bigl( X_i X_{i+|k|}, X_j X_{j+|\ell|} \bigr)\bigr|\\
&=& \frac{1}{N}\sum_{|r| < N} \biggl(1-\frac{|r|}{N}\biggr)
\bigl|\operatorname{Cov}\bigl( X_0 X_{|k|}, X_r X_{r+|\ell|} \bigr)\bigr|.
\end{eqnarray*}
The last sum can be split into three terms corresponding to
$r=0$, $r<0$ and $r>0$.

The contribution to the left-hand side of (\ref{e:lr3})
of the term corresponding to $r=0$ is
\[
N^{-1} \sum_{|k|, |\ell| \le q}
\bigl|\operatorname{Cov}\bigl( X_0 X_{|k|}, X_0 X_{|\ell|} \bigr)\bigr|= O(q^2/N).
\]
The terms corresponding to $r<0$ and $r>0$ are handled in the same way,
so we focus on the contribution of the summands with $r>0$
which is
\[
N^{-1} \sum_{|k|, |\ell| \le q} \sum_{r=1}^{N-1} \biggl(1-\frac
{r}{N}\biggr)
\bigl|\operatorname{Cov}\bigl( X_0 X_{|k|}, X_r X_{r+|\ell|} \bigr)\bigr|.
\]
We now use the decompositions
\begin{eqnarray*}
\operatorname{Cov}\bigl( X_0 X_{|k|}, X_r X_{r+|\ell|} \bigr)
&=& \operatorname{Cov}\bigl( X_0 X_{|k|}, X_r^{(r)} X_{r+|\ell|}^{(r+|\ell
|)} \bigr)\\
&&{} + \operatorname{Cov}\bigl( X_0 X_{|k|},
X_r X_{r+|\ell|} - X_r^{(r)} X_{r+|\ell|}^{(r+|\ell|)} \bigr)
\end{eqnarray*}
and
\begin{eqnarray*}
\operatorname{Cov}\bigl( X_0 X_{|k|}, X_r^{(r)} X_{r+|\ell|}^{(r+|\ell|)} \bigr)
&=& \operatorname{Cov}\bigl( X_0 X_{|k|}^{(|k|)}, X_r^{(r)} X_{r+|\ell
|}^{(r+|\ell|)}
\bigr)\\
&&{} + \operatorname{Cov}\bigl( X_0 \bigl(X_{|k|}- X_{|k|}^{(|k|)}\bigr),
X_r^{(r)} X_{r+|\ell|}^{(r+|\ell|)} \bigr).
\end{eqnarray*}
By Definition \ref{d:siegi}, $X_0$ depends on $\varepsilon_0,
\varepsilon_{-1},
\ldots$
while the random variables
$X_{|k|}^{(k)}, X_r^{(r)}$ and $X_{r+|\ell|}^{(r+|\ell|)}$ depend
on $\varepsilon_1, \varepsilon_2, \ldots, \varepsilon_{k \vee
(r+|\ell|)}$ and errors
independent of the $\varepsilon_i$. Therefore
$\operatorname{Cov}( X_0 X_{|k|}^{(|k|)}, X_r^{(r)} X_{r+|\ell
|}^{(r+|\ell|)}
)$ is equal to
\begin{eqnarray*}
&&E\bigl[ X_0 X_{|k|}^{(|k|)} X_r^{(r)} X_{r+|\ell|}^{(r+|\ell|)}\bigr]
- E \bigl[ X_0 X_{|k|} \bigr]
E \bigl[ X_r^{(r)} X_{r+|\ell|}^{(r+|\ell|)} \bigr]\\
&&\qquad =E [X_0] E\bigl[ X_{|k|}^{(|k|)} X_r^{(r)} X_{r+|\ell|}^{(r+|\ell
|)}\bigr]
- E [X_0] E \bigl[ X_{|k|}^{(|k|)} \bigr]
\bigl[ X_r^{(r)} X_{r+|\ell|}^{(r+|\ell|)} \bigr]= 0.
\end{eqnarray*}
We thus obtain
\begin{eqnarray*}
\operatorname{Cov}\bigl( X_0 X_{|k|}, X_r X_{r+|\ell|} \bigr)
&=& \operatorname{Cov}\bigl( X_0 \bigl(X_{|k|}- X_{|k|}^{(|k|)}\bigr),
X_r^{(r)} X_{r+|\ell|}^{(r+|\ell|)} \bigr)\\
&&{} + \operatorname{Cov}\bigl( X_0 X_{|k|},
X_r X_{r+|\ell|} - X_r^{(r)} X_{r+|\ell|}^{(r+|\ell|)} \bigr).
\end{eqnarray*}

By Assumption (\ref{e:lr-ct}), it remains
to verify that
\[
N^{-1} \sum_{|k|, |\ell| \le q} \sum_{r=1}^{N-1}
\bigl|\operatorname{Cov}\bigl( X_0 X_{|k|},
X_r X_{r+|\ell|} - X_r^{(r)} X_{r+|\ell|}^{(r+|\ell|)} \bigr)\bigr|\to0.
\]
This is done using the technique introduced in the
proof of Theorem \ref{t:ChatW}. By the Cauchy--Schwarz
inequality, the problem reduces to showing that
\[
N^{-1} \sum_{|k|, |\ell| \le q} \sum_{r=1}^{N-1}
\bigl\{ E\bigl[X_0^2 X_{|k|}^2\bigr] \bigr\}^{1/2}
\bigl\{ E \bigl[
\bigl( X_r X_{r+|\ell|} - X_r^{(r)} X_{r+|\ell|}^{(r+|\ell|)}\bigr)^2
\bigr]\bigr\}^{1/2} \to0.
\]
Using (\ref{e:abcd}), this in turn is bounded by constant times
\[
N^{-1} \sum_{|k|, |\ell| \le q} \sum_{r=1}^{\infty}
\bigl\{ E \bigl[ X_r - X_r^{(r)} \bigr]^4 \bigr\}^{1/4},
\]
which tends to zero by $L^4$--$m$-approximability and the condition
$q^2/N \to0$.
\end{pf*}
\begin{pf*}{Proof of Proposition \protect\ref{p:diffbartlett}}
We only
show the first part, the second is similar. Let $\omega_q(h)$ be the Bartlett
estimates satisfying Assumption \ref{a:wq}. Without loss of
generality we will assume below that the constant $b$ in
(\ref{e:og-b}) is 1.
Then the element in the $k$th row and
$\ell$th column of $\hat{{\bolds\Sigma}}(\bolds\beta)-\hat
{{\bolds\Sigma}}(\hat{\mathbf{C}
}\hat{\bolds\beta})$ is
\begin{eqnarray*}
&&\sum_{|h|\leq q}
\frac{\omega_q(h)}{N}\sum_{1\leq n\leq N-|h|}\bigl(\beta_{kn}\beta
_{\ell
,n+|h|}-{\hat c}_k
\hat{\beta}_{kn}{\hat c}_\ell\hat{\beta}_{\ell,n+|h|}\bigr)\\
&&\qquad =\sum_{|h|\leq q}
\frac{\omega_q(h)}{N}\sum_{1\leq n\leq N-|h|}\beta_{kn}\bigl(\beta
_{\ell,n+|h|}
-{\hat c}_\ell\hat{\beta}_{\ell,n+|h|}\bigr)\\
&&\qquad\quad{} + \sum_{|h|\leq q}
\frac{\omega_q(h)}{N}\sum_{1\leq n\leq N-|h|}{\hat c}_\ell\hat
{\beta}_{\ell
,n+|h|}(\beta_{kn}-{\hat c}_k\hat{\beta}_{kn})\\
&&\qquad=F_1(N,k,\ell)+F_2(N,k,\ell).
\end{eqnarray*}
For reasons of symmetry it is enough to estimate $F_1(N,k,\ell)$.
We have for any $t_N> 0$
\begin{eqnarray*}
&&P\bigl(|F_1(N,k,\ell)|>\varepsilon\bigr)\\
&&\qquad\leq\sum_{|h|\leq q}P \biggl(
\frac{\omega_q(h)}{N}\sum_{1\leq n\leq N-|h|}\beta_{kn}\bigl(\beta
_{\ell,n+|h|}
-{\hat c}_\ell\hat{\beta}_{\ell,n+|h|}\bigr)
>\frac{\varepsilon}{2q+1}
\biggr)\\
&&\qquad\leq\sum_{|h|\leq q}P \biggl(
\sum_{1\leq n\leq N-|h|}\beta_{kn}^2\sum_{1\leq n\leq N-|h|}\bigl(\beta
_{\ell,n+|h|}
-{\hat c}_\ell\hat{\beta}_{\ell,n+|h|}\bigr)^2
>\frac{\varepsilon^2 N^2}{(2q+1)^2} \biggr)\\
&&\qquad\leq(2q+1)P \biggl(
\sum_{1\leq n\leq N}\beta_{kn}^2>N(2q+1)t_N \biggr)\\
&&\qquad\quad{} +(2q+1)P \biggl(
\sum_{1\leq n\leq N}
(\beta_{\ell n}-{\hat c}_\ell\hat{\beta}_{\ell n})^2>\frac
{\varepsilon^2
N}{t_N(2q+1)^3}
\biggr)\\
&&\qquad=(2q+1)\bigl(P_1(k,N)+P_2(\ell,N)\bigr).
\end{eqnarray*}
By the Markov inequality and the fact that the $\beta_{kn}$, $1\leq
n\leq N$, are
identically distributed, we get for all $k\in\{1,\ldots,d\}$
\[
(2q+1)P_1(k,N)\leq\frac{E(\beta_{k1}^2)}{t_N}\leq\frac{E\|Y_1\|^2}{t_N},
\]
which tends to zero as long as $t_N\to\infty$.

The estimation of $P_2(\ell,N)$ requires a little bit more effort. We notice
first that
%
%e6.10 ###
%
\begin{equation}\label{e:varS2}
\limsup_{N\to\infty}
\frac{1}{N}\operatorname{Var} \biggl(\sum_{1\leq n\leq N}\|Y_n\|^2 \biggr)\leq
{\sum_{h\in\mathbb{Z}}}|{\operatorname{Cov}}(\|Y_1\|^2,\|Y_h\|^2)|<\infty.
\end{equation}
The summability of the latter series follows by now routine
estimates from (\ref{e:mL4}). For any $x,y>0$ we have
\begin{eqnarray*}
&&P \biggl(\sum_{1\leq n\leq N}
(\beta_{\ell n}-{\hat c}_\ell\hat{\beta}_{\ell n})^2>x
\biggr)\\
&&\qquad =P \biggl(
\sum_{1\leq n\leq N}
\biggl(\int Y_n(t)\bigl(v_\ell(t)-{\hat c}_\ell\hat{v}_\ell(t)\bigr)\,dt \biggr)^2>x
\biggr)\\
&&\qquad \leq P \biggl(
\sum_{1\leq n\leq N}
\|Y_n\|^2 \|v_\ell(t)-{\hat c}_\ell\hat{v}_\ell(t)\|^2>x
\biggr)\\
&&\qquad \leq P \biggl(
\sum_{1\leq n\leq N}
\|Y_n\|^2>xy \biggr)+P \bigl( \|v_\ell(t)-{\hat c}_\ell\hat{v}_\ell
(t)\|^2>x/y
\bigr)\\
&&\qquad=P_{21}(N)+P_{22}(\ell,N).
\end{eqnarray*}
If we require that
$y>NE\|Y_1\|^2/x$,
then by the Markov inequality and (\ref{e:varS2}) we have
\[
P_{21}(N)\leq\kappa_1 \biggl(\frac{xy}{\sqrt{N}}-\sqrt{N}E\|Y_1\|
^2 \biggr)^{-2}
\]
for some constant $\kappa_1$ which does not depend on $N$.
By Theorem \ref{t:B4W} and again the Markov inequality there exists a constant
$\kappa_2$ such that for all $\ell\in\{1,\ldots,d\}$
\[
P_{22}(\ell,N)\leq\kappa_2\frac{y}{xN}.
\]
The $x$ in the term $P_2(\ell,N)$ is given by
\[
x=\frac{\varepsilon^2N}{t_N(2q+1)^3}.
\]
Set $y=2NE\|Y_1\|^2/x$.
Then for all $\ell\in\{1,\ldots,d\}$
%
%e6.11 ###
%
\begin{equation}\label{e:P2}
P_{21}(N)\leq\kappa_1\frac{1}{(E\|Y_1\|^2)^2N} \quad\mbox{and}\quad
P_{22}(\ell,N)\leq\kappa_2\frac{2E\|Y_1\|^2}{\varepsilon
^4N^2}t_N^2(2q+1)^6.\hspace*{-32pt}
\end{equation}
Letting $t_N=(2q+1)^{1/2}$ shows that under $q^4/N\to0$ the term
$(2q+1)P_2(\ell,N)\hspace*{-0.41pt}\to0$. This finishes the proof of Proposition \ref
{p:diffbartlett}.
\end{pf*}

%s6.2 ###
\subsection{\texorpdfstring{Proofs of Theorems \protect\ref{t:QFFCLT} and
\protect\ref{th:alternative}}{Proofs of Theorems 5.1 and 5.2}} \label{s:pQ}

The proof of Theorem \ref{t:QFFCLT} relies on Theorem A.1 of
Aue et al. \cite{auehormannhorvathreimherr2009}, which we state here for
ease of reference.
\begin{theorem} \label{t:AA1}
Suppose $\{\bolds\xi_n\}$ is a $d$-dimensional $L^2$--$m$-approximable
mean zero sequence. Then
%
%e6.12 ###
%
\begin{equation}\label{e:AA1}
N^{-1/2} \mathbf{S}_N(\cdot, \bolds\xi) \stackrel{d}{\rightarrow
}\mathbf{W}({\bolds\xi})(\cdot),
\end{equation}
where $\{ \mathbf{W}({\bolds\xi})(x), x \in[0,1] \}$ is
a mean zero Gaussian process with covariances,
\[
\operatorname{Cov}( \mathbf{W}({\bolds\xi})(x), \mathbf{W}({\bolds
\xi})(y)) = \min(x,y) {\bolds\Sigma}(\bolds\xi).
\]
The convergence in (\ref{e:AA1}) is in the $d$-dimensional Skorokhod
space $D_d([0,1])$.
\end{theorem}
\begin{pf*}{Proof of Theorem \protect\ref{t:QFFCLT}}
Let
\[
G_N(x,\bolds\xi)=\frac{1}{N}\mathbf{L}_n(x,\bolds\xi)^T
\hat{{\bolds\Sigma}}(\bolds\xi)^{-1}\mathbf{L}_n(x,\bolds\xi)^T.
\]
We notice that
replacing the $\mathbf{L}_N(x,\hat{{\bolds\eta}})$ with $\mathbf
{L}_N(x,\hat{\bolds\beta})$ does not
change the test statistic in (\ref{e:test}).
Furthermore, since by the second part of Proposition \ref{p:diffbartlett}
$|\hat{{\bolds\Sigma}}(\hat{{\bolds\eta}})-\hat{{\bolds\Sigma
}}(\hat{\bolds\beta})|=o_P(1)$,
it is enough to study the limiting behavior of the
sequence $G_N(x,\hat{\bolds\beta})$.
This is done by first deriving the asymptotics
of $G_N(x,\bolds\beta)$ and then analyzing
the effect of replacing $\bolds\beta$ with $\hat{\bolds\beta}$.

Let $\bolds\beta_i^{(m)}$ be the $m$-dependent approximations for
$\bolds\beta_i$
which are obtained by replacing
$Y_i(t)$ in (\ref{e:beta}) by $Y_i^{(m)}(t)$. For a vector $\mathbf{v}$
in ${R}^d$ we let
$|\mathbf{v}|$ be its Euclidian norm. Then
\begin{eqnarray*}
E\bigl|\bolds\beta_1-\bolds\beta_1^{(m)}\bigr|^2
& = & E\sum_{\ell=1}^d \bigl(\beta_{\ell1}-\beta_{\ell1}^{(m)}
\bigr)^2\\
& = &\sum_{\ell=1}^d E \biggl(\int\bigl(Y_1(t)-Y_1^{(m)}(t)\bigr) v_\ell(t)
\,dt \biggr)^2\\
& \leq &\sum_{\ell=1}^dE\int\bigl(Y_1(t)-Y_1^{(m)}(t)\bigr)^2\,dt\int v_\ell
^2(t)\,dt\\
& = &d \nu_2^2\bigl(Y_1-Y_1^{(m)}\bigr).
\end{eqnarray*}
Since by Lyapunov's inequality we have
$\nu_2(Y_1-Y_1^{(m)})\leq\nu_4(Y_1-Y_1^{(m)})$,
(\ref{e:mL4}) yields that
$\sum_{m\geq1} (E|\bolds\beta_1-\bolds\beta
_1^{(m)}|^2)^{1/2}<\infty$.
Thus Theorem \ref{t:AA1} implies that
\[
\frac{1}{\sqrt{N}}\mathbf{S}_N(x,\bolds\beta)
\stackrel{D^d[0,1]}{\longrightarrow} \mathbf{W}(\bolds\beta)(x).
\]
The coordinatewise absolute convergence of the series
${\bolds\Sigma}(\bolds\beta)$ follows from part (a) of
Theorem \ref{t:lr-v}. By assumption the estimator
$\hat{{\bolds\Sigma}}(\bolds\beta)$ is consistent, and
consequently
\[
\int G_N(x,\bolds\beta)\,dx\stackrel{D[0,1]}{\longrightarrow}\sum
_{\ell=1}^d
\int B_\ell^2(x)\,dx
\]
follows from the continuous mapping theorem.

We turn now to the effect of changing $G_N(x,\bolds\beta)$
to $G_N(x,\hat{\bolds\beta})$.
Due to the quadratic structure of
$G_N(x,\bolds\xi)$, we have $G_N(x,\hat{\bolds\beta})=G_N(x,\hat
{\mathbf{C}}\hat{\bolds\beta})$ when
$\hat{\mathbf{C}}=\operatorname{diag}({\hat c}_1,{\hat c}_2,\ldots
,{\hat c}_d)$.
To finish the proof it is thus sufficient to show that
%
%e6.13 ###
%
\begin{equation}\label{diff1}
\sup_{x\in[0,1]}\frac{1}{\sqrt{N}}
|\mathbf{S}_N(x,\bolds\beta)-\mathbf{S}_N(x,\hat{\mathbf{C}}\hat
{\bolds\beta})|=o_P(1)
\end{equation}
and
%
%e6.14 ###
%
\begin{equation}\label{diff2}
|\hat{{\bolds\Sigma}}(\bolds\beta)-\hat{{\bolds\Sigma}}(\hat
{\mathbf{C}}\hat{\bolds\beta})|=o_P(1).
\end{equation}
Relation (\ref{diff2}) follows from Proposition \ref{p:diffbartlett}.
To show (\ref{diff1}) we observe that
by the Cauchy--Schwarz inequality and Theorem \ref{t:B4W}
\begin{eqnarray*}
&&\sup_{x\in[0,1]}\frac{1}{N}
|\mathbf{S}_N(x,\bolds\beta)-\mathbf{S}_N(x,\hat{\mathbf{C}}\hat
{\bolds\beta})|^2\\
&&\qquad =\sup_{x\in[0,1]}\frac{1}{N}
\sum_{\ell=1}^d \Biggl|\int
\sum_{n=1}^{\lfloor Nx\rfloor}
Y_n(t)\bigl(v_\ell(t)-{\hat c}_\ell{\hat v}_\ell(t)\bigr)\,dt \Biggr|^2\\
&&\qquad \leq\frac{1}{N}\sup_{x\in[0,1]}\int\Biggl(
\sum_{n=1}^{\lfloor Nx\rfloor}
Y_n(t) \Biggr)^2 \,dt\times\sum_{\ell=1}^d\int
\bigl(v_\ell(t)-{\hat c}_\ell{\hat v}_\ell(t)\bigr)^2\,dt\\
&&\qquad \leq\frac{1}{N}\int\max_{1\leq k\leq N} \Biggl(
\sum_{n=1}^{k}
Y_n(t) \Biggr)^2 \,dt\times O_P(N^{-1}).
\end{eqnarray*}
Define
\[
g(t)=E|Y_1(t)|^2+
2 (E|Y_1(t)|^2 )^{1/2}\sum_{r\geq1}
\bigl(E\bigl|Y_{1+r}(t)-Y_{1+r}^{(r)}(t)\bigr|^2
\bigr)^{1/2}.
\]
Then by similar arguments as in Section \ref{s:p-si} we have
\[
E \Biggl(
\sum_{n=1}^{N}
Y_n(t) \Biggr)^2\leq N g(t).
\]
Hence by Menshov's inequality (see, e.g., Billingsley \cite
{billingsley1999}, Section 10) we infer that
\[
E\max_{1\leq k\leq N} \Biggl(
\sum_{n=1}^{k}
Y_n(t) \Biggr)^2\leq(\log\log4N)^2 N g(t).
\]
Notice that (\ref{e:mL4}) implies $\int g(t)\,dt<\infty$. In turn we
obtain that
\[
\frac{1}{N}\int\max_{1\leq k\leq N} \Biggl(
\sum_{n=1}^{k}
Y_n(t) \Biggr)^2 \,dt=O_P ((\log\log N)^2 ),
\]
which proves (\ref{diff1}).
\end{pf*}
\begin{pf*}{Proof of Theorem \protect\ref{th:alternative}}
Notice that if the mean function changes from $\mu_1(t)$ to $\mu_2(t)$
at time $k^*=\lfloor N\theta\rfloor$,
then $\mathbf{L}_N(x,\bolds{\hat\eta})$ can be written as
%
%e6.15 ###
%
\begin{equation}\label{e:cpSN}
\mathbf{L}_N(x,\bolds{\hat\beta})+
N\cases{
x(1-\theta)[\hat{\bolds\mu}_1-\hat{\bolds\mu}_2], &\quad if
$x\leq\theta$;\cr
\theta(1-x)[\hat{\bolds\mu}_1-\hat{\bolds\mu}_2], &\quad if
$x> \theta$,}
\end{equation}
where
\[
\hat{\bolds\mu}_1 = \biggl[\int\mu_1(t) {\hat v}_1(t)\,dt, \int\mu_1(t)
{\hat v}_2(t)\,dt,
\ldots,\int\mu_1(t) {\hat v}_d(t)\,dt\biggr]^T
\]
and $\hat{\bolds\mu}_2$ is defined analogously.

It follows from (\ref{e:cpSN}) that $T_N(d)$ can be
expressed as the sum of three terms:
\[
T_N(d) = T_{1,N}(d) + T_{2,N}(d) + T_{3,N}(d),
\]
where
\begin{eqnarray*}
T_{1,N}(d)
&=& \frac{1}{N}\int_0^1\mathbf{L}_N(x,\bolds{\hat\beta})^T
\hat{\bolds\Sigma}(\hat{\bolds\eta})^{-1}
\mathbf{L}_N(x,\bolds{\hat\beta})
\,dx;
\\
T_{2,N}(d) &=& \frac{N}{2} \theta(1-\theta)[\hat{\bolds\mu
}_1-\hat{\bolds\mu}_2]^T
\hat{\bolds\Sigma}(\hat{\bolds\eta})^{-1}[\hat{\bolds
\mu}
_1-\hat{\bolds\mu}_2];
\\
T_{3,N}(d) &=& \int_0^1 g(x, \theta)
\mathbf{L}_N(x,\bolds{\hat\beta})^T
\hat{\bolds\Sigma}(\hat{\bolds\eta})^{-1}[\hat{\bolds
\mu}
_1-\hat{\bolds\mu}_2]\,dx,
\end{eqnarray*}
with
$
g(x, \theta) = 2 \{x(1-\theta)I_{\{x\leq\theta\}}
+\theta(1-x)I_{\{x> \theta\}} \}.
$

Since $\bolds\Omega$ in (\ref{e:conalt})
is positive definite (p.d.),
$\hat{\bolds\Sigma}(\hat{\bolds\eta})$ is almost surely
p.d. for large enough $N$ ($N$ is random).
Hence for large enough $N$ the term $T_{1,N}(d)$ is nonnegative. We
will show that
$N^{-1}T_{2,N}(d) \ge\kappa_1 + o_P(1)$, for a positive constant
$\kappa_1$,
and $N^{-1}T_{3,N}(d) = o_P(1)$. To this end we notice the following.
Ultimately all eigenvalues of $\hat{\bolds\Sigma}(\hat
{\bolds\eta})$
are positive. Let $\lambda^*(N)$ and $\lambda_*(N)$ denote the largest,
respectively, the smallest eigenvalue.
By Lemma \ref{l:B4.2},
$\lambda^*(N)\to\lambda^*$ a.s. and $\lambda_*(N)\to\lambda_*$
a.s., where $\lambda^*$ and $\lambda_*$ are the largest
and smallest eigenvalue of
$\bolds\Omega$. Next we claim that
\[
|\hat{\bolds\mu}_1-\hat{\bolds\mu}_2| = |\bolds\mu_1-\bolds\mu
_2| + o_P(1).
\]
To obtain this, we use the relation $\|{\hat v}_i - {\hat c}_j v_j\|
= o_P(1)$ which can be proven similarly as Lemma A.1 of
Berkes et al. \cite{berkesgabryshorvathkokoszka2009}, but
the law of large numbers in a Hilbert space must be replaced
by the ergodic theorem.
The ergodicity of $\{ Y_n\}$ follows
from the representation $Y_n = f(\varepsilon_n, \varepsilon_{n-1},
\ldots)$.
Notice that because of the presence of
a change point it cannot be claimed that
$\|{\hat v}_i - {\hat c}_j v_j\| = O_P(N^{-1/2})$.

It follows that if $N$ is large enough, then
\[
[\hat{\bolds\mu}_1-\hat{\bolds\mu}_2]^T\hat{\bolds\Sigma}
(\hat{\bolds\eta})^{-1}[\hat{\bolds\mu}_1-\hat{\bolds\mu
}_2]>\frac{1}{2\lambda^*}
|\hat{\bolds\mu}_1-\hat{\bolds\mu}_2|^2 = \frac{1}{2\lambda
^*}|\bolds\mu_1-\bolds\mu_2|^2 + o_P(1).
\]

To verify $N^{-1}T_{3,N}(d) = o_P(1)$, observe that
\begin{eqnarray*}
&&{\sup_{x\in[0,1]}} |\mathbf{L}_N(x,\bolds{\hat\beta})^T
\hat{\bolds\Sigma}(\hat{\bolds\eta})^{-1}[\hat{\bolds
\mu}
_1-\hat{\bolds\mu}_2] |\\
&&\qquad \leq{\sup_{x\in[0,1]}}|\mathbf{L}_N(x,\bolds{\hat\beta
})| \times
|\hat{\bolds\Sigma}(\hat{\bolds\eta})^{-1}| \times
|\hat{\bolds\mu}_1-\hat{\bolds\mu}_2|\\
&&\qquad = o_P(N) |\bolds\mu_1-\bolds\mu_2|.
\end{eqnarray*}
We used the matrix norm $|A|={\sup_{|x|\leq1}}|Ax|$ and
$|\hat{\bolds\Sigma}(\hat{\bolds\eta})^{-1}|
\stackrel{\mathrm{a.s}}{\longrightarrow}|\bolds{\Omega
}^{-1}|<\infty$.
\end{pf*}

%s6.3 ###
\subsection{\texorpdfstring{Proof of Theorem \protect\ref{t:linregdep}}{Proof of Theorem 6.1}}\label{s:p-linreg}
We first establish a technical bound which implies
the consistency of the estimator
$\hat\sigma_{\ell k}$ given in (\ref{e:hat-sigma-lk}).
Let $\hat{c}_\ell=\operatorname{sign}(\langle v_\ell,\hat{v}_\ell
\rangle)$ and
$\hat{d}_k=\operatorname{sign}(\langle u_k,\hat{u}_k \rangle)$.
\begin{lemma}\label{l:consistentsigma}
Under the assumptions of Theorem \ref{t:linregdep}
we have
\[
\limsup_{N\to\infty} NE|\sigma_{\ell k}-\hat{c}_\ell \hat
{d}_k \hat\sigma_{\ell k}|^2\leq
\kappa_1 \biggl(\frac{1}{\alpha_k^2}+\frac{1}{(\alpha_\ell
')^2} \biggr),
\]
where $\kappa_1$ is a constant independent of $k$ and $\ell$.
\end{lemma}
\begin{pf}
It follow from elementary inequalities that
\[
|\sigma_{\ell k}-\hat{c}_\ell \hat{d}_k \hat\sigma_{\ell
k}|^2\leq2T_1^2+2T_2^2,
\]
where
\begin{eqnarray*}
T_1&=&\frac{1}{N}\dint\Biggl(\sum_{i=1}^N
\bigl(X_i(s)Y_i(t)-E[X_i(s)Y_i(t)] \bigr) \Biggr)u_k(s)v_\ell(t)\,dt\,ds;\\
T_2&=&\frac{1}{N}\sum_{i=1}^N\dint E[X_i(s)Y_i(t)]
[u_k(t)v_\ell(s)-\hat{d}_k\hat{u}_k(t)
\hat{c}_\ell\hat{v}_\ell(s) ]\,dt\,ds.
\end{eqnarray*}
By the Cauchy--Schwarz inequality and (\ref{e:abcd}) we obtain
\begin{eqnarray*}
T_1^2&\leq&\frac{1}{N^2}\dint\Biggl(\sum_{i=1}^N
X_i(s)Y_i(t)-E[X_i(s)Y_i(t)] \Biggr)^2\,dt\,ds;\\
T_2^2 &=& 2\nu_2^2(X)\nu_2^2(Y) (\|u_k-\hat{d}_k\hat{u}_k\|^2+\|
v_\ell-\hat{c}_\ell\hat{v}_\ell\|^2 ).
\end{eqnarray*}
Hence by similar arguments as we used for the proof of Theorem \ref{t:ChatW}
we get $NET_1^2=O(1)$.
The proof follows now immediately from Lemma \ref{l:B4.3} and
Theorem \ref{t:ChatW}.
\end{pf}

Now we are ready to verify (\ref{e:ymw-con}).
We have
\[
\hat\psi_{KL}(t,s)=\sum_{k=1}^K\sum_{\ell=1}^L\hat\lambda_\ell
^{-1}\hat\sigma_{\ell k}
\hat{u}_k(t)\hat{v}_\ell(s).
\]
The orthogonality of the sequences $\{u_k\}$ and $\{v_\ell\}$
and (\ref{e:ymw}) imply that
\begin{eqnarray*}
&&\dint \biggl(
\sum_{k>K}\sum_{\ell>L}\lambda_\ell^{-1}\sigma_{\ell k}
u_k(t)v_\ell(s)
\biggr)^2\,dt\,ds\\
&&\qquad=\sum_{k>K}\sum_{\ell>L}\dint\lambda_\ell^{-2}\sigma_{\ell k}^2
u_k^2(t)v_\ell^2(s)\,dt\,ds\\
&&\qquad=\sum_{k>K}\sum_{\ell>L}\lambda_\ell^{-2}\sigma_{\ell k}^2\to
0\qquad (L,K\to\infty).
\end{eqnarray*}
Therefore, letting
\[
\psi_{KL}(t,s)=\sum_{k=1}^K\sum_{\ell=1}^L\lambda_\ell^{-1}\sigma
_{\ell k}u_k(t)v_\ell(s),
\]
(\ref{e:ymw-con}) will follow once we show that
\[
\dint[\psi_{KL}(t,s)-\hat\psi_{KL}(t,s)
]^2\,dt\,ds\stackrel{P}{\rightarrow}
0\qquad (N\to\infty).
\]
Notice that by the Cauchy--Schwarz inequality the latter relation is
implied by
%
%e6.16 ###
%
\begin{eqnarray}\label{e:toshow}\qquad
KL\sum_{k=1}^K\sum_{\ell=1}^L\dint
[\lambda_\ell^{-1}\sigma_{\ell k}u_k(t)v_\ell(s)-
\hat\lambda_\ell^{-1}\hat\sigma_{\ell k}\hat{u}_k(t)\hat{v}_\ell
(s) ]^2\,dt\,ds
\stackrel{P}{\rightarrow}0\nonumber\\[-8pt]\\[-8pt]
\eqntext{(N\to\infty).}
\end{eqnarray}
A repeated application of (\ref{e:abcd}) and some basic algebra yield
\begin{eqnarray*}
&&\tfrac{1}{4} [\lambda_\ell^{-1}\sigma_{\ell k}u_k(t)v_\ell(s)-
\hat\lambda_\ell^{-1}\hat\sigma_{\ell k}\hat{u}_k(t)\hat{v}_\ell
(s) ]^2\\
&&\qquad \leq
\lambda_\ell^{-2}|\sigma_{\ell k}-
\hat{c}_\ell \hat{d}_k \hat\sigma_{\ell k}|^2\hat
{u}_k^2(t)\hat{v}_\ell^2(s)
+\hat\sigma_{\ell k}^2 |\lambda_\ell^{-1}-\hat\lambda_\ell
^{-1} |^2
\hat{u}_k^2(t)\hat{v}_\ell^2(s)\\
&&\qquad\quad{} +\sigma_{\ell k}^2\lambda_\ell^{-2}
|u_k(t)-\hat{d}_k\hat{u}_k(t) |^2v_\ell^2(s)
+\sigma_{\ell k}^2\lambda_\ell^{-2} |v_\ell(s)-\hat{c}_\ell
\hat{v}_\ell(s) |^2
\hat{u}_k^2(t).
\end{eqnarray*}
Hence
\begin{eqnarray*}
&&\tfrac{1}{4}\dint[\lambda_\ell^{-1}\sigma_{\ell
k}u_k(t)v_\ell(s)-
\hat\lambda_\ell^{-1}\hat\sigma_{\ell k}\hat{u}_k(t)\hat{v}_\ell
(s) ]^2\,dt\,ds\\
&&\qquad \leq
\lambda_\ell^{-2}|\sigma_{\ell k}-
\hat{c}_\ell \hat{d}_k \hat\sigma_{\ell k}|^2
+\hat\sigma_{\ell k}^2 |\lambda_\ell^{-1}-\hat\lambda_\ell
^{-1} |^2
\\
&&\qquad\quad{} +\sigma_{\ell k}^2\lambda_\ell^{-2} (
\|u_k-\hat{d}_k\hat{u}_k\|^2
+\|v_\ell-\hat{c}_\ell\hat{v}_\ell\|^2 ).
\end{eqnarray*}
Thus in order to get (\ref{e:toshow}) we will show that
%
%e6.19 ###
%e6.18 ###
%e6.17 ###
%
\begin{eqnarray}
\label{e:ts1}
KL\sum_{k=1}^K\sum_{\ell=1}^L\lambda_\ell^{-2}|\sigma_{\ell k}-
\hat{c}_\ell \hat{d}_k\hat\sigma_{\ell k}|^2&\stackrel
{P}{\rightarrow}&0;\\
\label{e:ts2}
KL\sum_{k=1}^K\sum_{\ell=1}^L\hat\sigma_{\ell k}^2 |\lambda
_\ell^{-1}-\hat\lambda_\ell^{-1} |^2
&\stackrel{P}{\rightarrow}& 0;\\
\label{e:ts3}
KL\sum_{k=1}^K\sum_{\ell=1}^L\sigma_{\ell k}^2\lambda_\ell
^{-2} (\|u_k-\hat{d}_k\hat{u}_k\|^2
+\|v_\ell-\hat{c}_\ell\hat{v}_\ell\|^2 )&\stackrel
{P}{\rightarrow}&0.
\end{eqnarray}
We start with (\ref{e:ts1}). By Lemma \ref{l:consistentsigma} and
Assumption \ref{a:linreg} we
have
\[
E \Biggl(KL\sum_{k=1}^K\sum_{\ell=1}^L\lambda_\ell^{-2}|
\sigma_{\ell k}-\hat{c}_\ell \hat{d}_k\hat{\sigma}_{\ell
k}|^2 \Biggr)\to0\qquad (N\to\infty).
\]

Next we prove relation (\ref{e:ts2}). In order to shorten the proof
we replace $\hat\sigma_{\ell k}$ by $\sigma_{\ell k}$. Otherwise we
would need a further intermediate step, requiring similar
arguments which follow. Now for any $0<\varepsilon<1$ we have
\begin{eqnarray*}
&& P\Biggl(KL\sum_{k=1}^K\sum_{\ell=1}^L
\sigma_{\ell k}^2 |\lambda_\ell^{-1}-\hat\lambda_\ell
^{-1} |^2>\varepsilon
\Biggr)
\\
&&\qquad=
P \Biggl(KL\sum_{k=1}^K\sum_{\ell=1}^L
\sigma_{\ell k}^2\lambda_\ell^{-2} \biggl|\frac{\hat\lambda_\ell
-\lambda_\ell}
{\hat\lambda_\ell} \biggr|^2>\varepsilon
\Biggr)\\
&&\qquad\leq
P \biggl(\max_{1\leq\ell\leq L}
\biggl|\frac{\hat\lambda_\ell-\lambda_\ell}
{\hat\lambda_\ell} \biggr|^2>\frac{\varepsilon}{\Psi KL}
\biggr)
\\
&&\qquad\leq
\sum_{\ell=1}^L
P \biggl(
\biggl|\frac{\hat\lambda_\ell-\lambda_\ell}
{\hat\lambda_\ell} \biggr|^2>\frac{\varepsilon}{\Psi KL}\cap
|\lambda_\ell-\hat\lambda_\ell|<\varepsilon\lambda_\ell
\biggr)
\\
&&\qquad\quad{} +\sum_{\ell=1}^L
P \biggl(
\biggl|\frac{\hat\lambda_\ell-\lambda_\ell}
{\hat\lambda_\ell} \biggr|^2>\frac{\varepsilon}{\Psi KL}\cap
|\lambda_\ell-\hat\lambda_\ell|\geq\varepsilon\lambda_\ell
\biggr)
\\
&&\qquad\leq
\sum_{\ell=1}^L \biggl[
P \biggl(
|\hat\lambda_\ell-\lambda_\ell|^2>\frac{\varepsilon}{\Psi KL}
\lambda_\ell(1-\varepsilon) \biggr)+
P (|\lambda_\ell-\hat\lambda_\ell|^2\geq\varepsilon^2\lambda
_\ell^2 )
\biggr]
\\
&&\qquad\leq
\kappa_2 \biggl(\frac{KL^2}{\epsilon N\lambda_L}+\frac
{1}{\varepsilon N\lambda_L^2} \biggr),
\end{eqnarray*}
by an application of the Markov inequality and Theorem \ref{t:B4W}.
According to our Assumption \ref{a:linreg} this also goes to zero for
$N\to\infty$.

Finally we prove (\ref{e:ts3}).
By Lemma \ref{l:B4.3} and Theorem \ref{t:ChatW}
we infer that
\begin{eqnarray*}
&&E \Biggl(KL\sum_{k=1}^K\sum_{\ell=1}^L\sigma_{\ell k}^2\lambda
_\ell^{-2} (
\|u_k-\hat{d}_k\hat{u}_k\|^2
+\|v_\ell-\hat{c}_\ell\hat{v}_\ell\|^2 ) \Biggr)\\
&&\qquad\leq\kappa
_3\frac{KL}{N}
\sum_{k=1}^K\sum_{\ell=1}^L\sigma_{\ell k}^2\lambda_\ell
^{-2} \biggl(\frac{1}{\alpha_k^2}
+\frac{1}{{\alpha_\ell'}^2} \biggr)\\
&&\qquad\leq2\kappa_3\Psi\frac{KL}{N\min\{h_L,h_K'\}^2}.
\end{eqnarray*}
Assumption \ref{a:linreg}(ii) assures that the last term goes to zero.
The proof is now complete.
\end{appendix}

\section*{Acknowledgments}
We thank P. Hall, H.-G. M{\"u}ller, D. Paul and J.-L. Wang for their
useful comments. We are grateful to two referees for raising a number
of points which helped us clarify several important issues.
Finally, we thank the Associate Editor for constructive criticism and
clear guidelines.

\printaddresses

\end{document}